\documentclass[journal]{IEEEtran}

\usepackage{graphicx}
\usepackage{color}
\usepackage{amsmath}
\usepackage{amssymb}
\usepackage{mathrsfs}
\usepackage{dsfont}
\usepackage{cite}  
\usepackage{epstopdf}   		
\usepackage{booktabs}
\usepackage[colorlinks=true,linkcolor=black,citecolor=black]{hyperref}

\usepackage[font=footnotesize]{caption}
\usepackage{subcaption}

\graphicspath{{./Figures/}}

\newtheorem{remark}{Remark}

\newcommand\oprocendsymbol{\hbox{$\square$}}
\newcommand\oprocend{\relax\ifmmode\else\unskip\hfill\fi\oprocendsymbol}

\begin{document}
%

\title{Secondary Frequency and Voltage Control of Islanded Microgrids via Distributed Averaging}

\author{
John~W.~Simpson-Porco,~\IEEEmembership{Student Member,~IEEE},
Qobad~Shafiee,~\IEEEmembership{Member,~IEEE},
Florian~D\"{o}rfler,~\IEEEmembership{Member,~IEEE},
Juan~C.~Vasquez,~\IEEEmembership{Senior Member,~IEEE},
Josep~M.~Guerrero,~\IEEEmembership{Fellow,~IEEE}, 
and~Francesco~Bullo,~\IEEEmembership{Fellow,~IEEE}

\thanks{Manuscript received October 14, 2014; revised January 22, 2015 and March 20, 2015; accepted April 24, 2015.} 
\thanks{Copyright (c) 2015 IEEE. Personal use of this material is permitted. However, permission to use this material for any other purposes must be obtained from the IEEE by sending a request to pubs-permissions@ieee.org.}
\thanks{This work was supported by NSF grants CNS-1135819 and by the National Science and Engineering Research Council of Canada.}
\thanks{J.~W.~Simpson-Porco and F.~Bullo are with the Center for Control, Dynamical Systems and Computation, University of California at Santa Barbara. Email: {\{johnwsimpsonporco,} {bullo\}@engineering.ucsb.edu}.}
\thanks{Q.~Shafiee, J.~C.~Vasquez, and J.~M.~Guerrero are with the Institute of Energy Technology, Aalborg University, Aalborg East DK-9220, Denmark. Email: {\{qsh, juq, joz\}@et.aau.dk}.}
\thanks{F. D\"orfler is with the Automatic Control Laboratory, Swiss Federal Institute (ETH) Zurich 
{dorfler@ethz.ch}.}
}

\markboth{IEEE Transactions on Industrial Electronics}%
{Shell \MakeLowercase{\textit{et al.}}: Bare Demo of IEEEtran.cls for Journals}
%

\maketitle


\begin{abstract}
In this work we present new distributed controllers for secondary frequency and voltage control in islanded microgrids. 
Inspired by techniques from cooperative control, the proposed controllers use localized information and nearest-neighbor communication to collectively perform secondary control actions.
The frequency controller rapidly regulates the microgrid frequency to its nominal value while maintaining active power sharing among the distributed generators. Tuning of the voltage controller provides a simple and intuitive trade-off between the conflicting goals of voltage regulation and reactive power sharing.
{ Our designs require no knowledge of the microgrid topology, impedances or loads.} 
The distributed architecture allows for flexibility and redundancy, and eliminates the need for a central microgrid controller. 
{ We provide a voltage stability analysis and} present extensive experimental results validating our designs, verifying robust performance under communication failure and during plug-and-play operation.
\end{abstract}


%
%
%




%
\IEEEpeerreviewmaketitle

\section{Introduction}
\label{Section: Introduction}


%
%
\IEEEPARstart{E}{}conomic factors, environmental concerns, and the rapidly expanding integration of small-scale renewable energy sources are all pushing the incumbent centralized power generation paradigm towards a more distributed future.
As a bridge between high-voltage transmission and low-voltage distributed generation (DG), the concept of a \emph{microgrid} continues to gain popularity \cite{RHL:02,JMG-JCV-JM-LGDV-MC:11,QCZ-TH:13,JMG-MC-TL-PCL:13a,JMG-MC-TL-PCL:13b}. Microgrids are low-voltage electrical distribution networks, heterogeneously composed of distributed generation, storage, load, and managed autonomously from the larger primary network. Microgrids can connect to a larger power system through a Point of Common Coupling (PCC), but are also able to ``island'' themselves and operate independently \cite{JMG-JCV-JM-LGDV-MC:11}. Islanded operation of a microgrid could be planned, or could occur spontaneously if a fault triggers the disconnection of the microgrid from the primary grid.

\begin{figure}[t]
\begin{center}
\includegraphics[width=0.7\columnwidth]{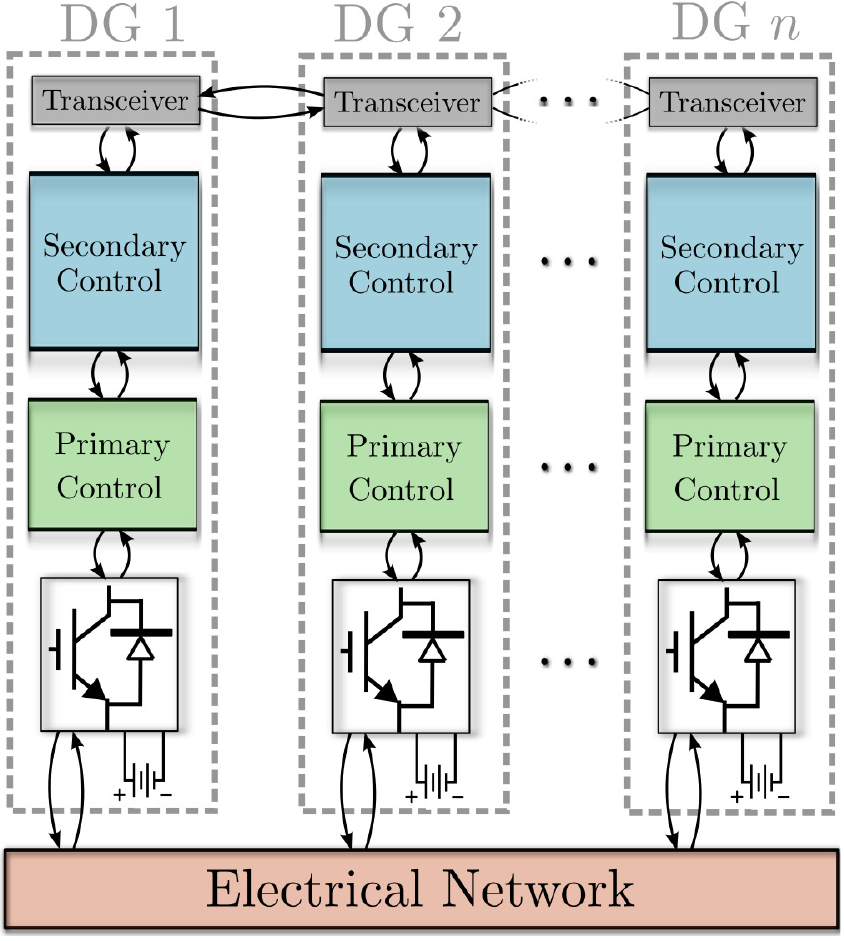}
\captionsetup{justification=raggedright,singlelinecheck=false}
\caption{Low-detail schematic of microgrid control architecture.
}
\label{Fig:DAPIFeedbackSimple}
\end{center}  
\end{figure}

Energy generation within a microgrid can be quite heterogeneous, including photovoltaic, wind, micro-turbines, etc. Such sources generate either DC power or variable frequency AC power, and are interfaced with a synchronous AC grid via power electronic \emph{inverters}. It is through these inverters that cooperative actions must be taken to ensure synchronization, voltage regulation, power balance and load sharing in the network \cite{JAPL-CLM-AGM:06}. Control strategies ranging from centralized to completely decentralized have been proposed to address these challenges \cite{AM-EO-DM-OO:10}, and have subsequently been aggregated into a hierarchical control architecture \cite{JMG-JCV-JM-LGDV-MC:11} (Figure \ref{Fig:DAPIFeedbackSimple}).

The control hierarchy consists of three levels. The first and most basic level is \emph{primary} control, which stabilizes the microgrid and establishes power sharing. Although centralized architectures have been used for primary control \cite{AM-EO-DM-OO:10}, in order to enhance redundancy and enable plug-and-play functionality, the current standard is to employ proportional control loops locally at each inverter. While successful for stabilization, these decentralized ``droop'' controllers force the bus voltages and the steady-state network frequency to deviate from their nominal values \cite{MCC-DMD-RA:93,JAPL-CLM-AGM:06,JMG-JCV-JM-MC-LGDV:09}. 

This leads naturally to the next level in the hierarchy, termed \emph{secondary} control. Broadly speaking, the goal of secondary control is to remove the aforementioned deviations in both global frequency and local voltage \cite{JAPL-CLM-AGM:06}. Centralized techniques for secondary control have been well studied in high-voltage transmission and distribution networks \cite{JM-JWB-JRB:08}. These centralized strategies have also been applied in the context of microgrids, and the term ``secondary'' has been broadened to include additional control goals such as reactive power sharing \cite{AM-MA-CSS-JMG-JCV:14,QS-JMG-JCV:14}, harmonic compensation, and voltage unbalance \cite{JMG-JCV-JM-LGDV-MC:11,JAPL-CLM-AGM:06,MCC-DMD-RA:93,MS-AJ-JCV-JGM:12b}. The final level of \emph{tertiary} control is concerned with global economic dispatch over the network, and depends on current energy markets and prices.
%
%

{ Several recent works (see Section \ref{Subsection: Review of Secondary Control}) have proposed secondary control strategies for microgrids. However, to date no single control strategy has been able to offer a flexible, plug-and-play architecture guaranteeing frequency and voltage regulation while maintaining precise active and reactive power sharing among non-identical DGs. Currently, this combination of goals appears infeasible with decentralized control using only local information (voltage, power, ect.) at each DG \cite{MA-DVD-HS-KHJ:14b}. As such, \emph{communication between DGs} has been identified as a key ingredient in achieving these
 goals while avoiding a centralized control architecture \cite{JWSP-FD-FB:12u,JWSP-FD-QS-JMG-FB:13e,HSVSKN-SD:13,MA-DVD-HS-KHJ:14b,QCZ:13,HB-JWSP-FD-FB:13j,QS-JMG-JCV:14}.}


{ In this paper we build on our previous theoretical and experimental works \cite{JWSP-FD-FB:12u,JWSP-FD-QS-JMG-FB:13e} and introduce a general and fully distributed framework for secondary frequency and voltage control in islanded microgrids. 
Our designs overcome the limitations of existing strategies by combining decentralized proportional droop control and integral control with distributed averaging algorithms. We therefore refer to our proposed controllers as \textbf{DAPI} (Distributed Averaging Proportional Integral) controllers. 
These controllers use decentralized control actions and sparse communication among neighboring DG units to achieve precise frequency regulation, active power sharing, and a tunable trade-off between voltage regulation and reactive power sharing. 
The distributed architecture eliminates the need for a central supervisory control: additional DGs are integrated through a low-bandwidth communication link to an existing DG, with the communication topology being a tunable design variable. The DAPI controllers are model-free, in the sense that they require no \emph{a priori} knowledge of the microgrid topology, line impedances or load demands. 

There are four main technical contributions in this paper. 
First, In Section \ref{Sec:ReactiveSharing} we highlight and clearly demonstrate a fundamental limitation of voltage control: precise voltage regulation and precise reactive power sharing are conflicting objectives. The presentation frames and motivates our subsequent controller designs. 
Second, in Section \ref{Section:DAPI} we review the frequency DAPI controller \cite{JWSP-FD-FB:12u} and introduce the voltage DAPI controller. This new voltage DAPI controller accounts for the conflict between voltage regulation and reactive power sharing by allowing for a tunable compromise between the two objectives. We build intuition for our design by detailing several tuning strategies. 
Taken together, the two DAPI controllers form a distributed duo for plug-and-play microgrid control.
{Third, in Section \ref{Sec:Stability} we present a small-signal voltage stability analysis of the microgrid under DAPI control, derive sufficient conditions on the controller gains and microgrid parameters for closed-loop stability, and study the transient performance of the system under changes in the controller gains.}
Finally, in Section \ref{Sec:Experiments} we present extensive experimental results validating our DAPI designs. The experimental microgrid consists of four heterogeneous DGs in a non-parallel configuration, with high $R/X$ connections and distributed load. We validate our designs, and move beyond our theoretical results by demonstrating controller performance under communication link failures and plug-and-play operation.

Section \ref{Section: Review of Distributed Control in Microgrids} contains a review of standard control strategies for microgrids. In particular, Section \ref{Subsection: Review of Primary Droop Controller} reviews primary droop control, while Section \ref{Subsection: Review of Secondary Control} provides a detailed review of secondary control for both frequency and voltage. Our main contributions are housed in Sections \ref{Sec:ReactiveSharing}--\ref{Sec:Experiments}, with concluding remarks being offered in Section \ref{Sec:Conclusions}.}
%
\section{Review of Microgrid Control}
\label{Section: Review of Distributed Control in Microgrids}
%
\subsection{Problem Setup and Review of Power Flow}
\label{Subsection: Review of Conventional Droop Controller}

In this work we consider microgrids consisting of $n$ buses which are either DGs or loads. For inductive lines of reactance $X_{ij}$ connecting bus $i$ to bus $j$, the active and reactive power injections $P_i$ and $Q_i$ at bus $i$ are given by \cite{JM-JWB-JRB:08}
\begin{subequations}\label{eq: power flow}
\begin{align}
	P_i &= \sum_{j=1}^n \nolimits \frac{E_iE_j}{X_{ij}}\sin(\theta_i-\theta_j)
	\,,
	\label{eq: power flow -- active}\\
	\Bigl.
	Q_i &= \frac{E_i^2}{X_i}-\sum_{j=1}^n \nolimits \frac{E_iE_j}{X_{ij}}\cos(\theta_i-\theta_j)
	\,,
	\label{eq: power flow -- reactive}
\end{align}
\end{subequations}
where $E_i$ (resp. $\theta_i$) is the voltage magnitude (resp. voltage phase angle) at bus $i$ and $X_{i} = 1/(\sum_{j=1}^n \nolimits X_{ij}^{-1})$.
\subsection{Review of Primary Droop Control}
\label{Subsection: Review of Primary Droop Controller}

A complete survey of primary control is beyond the scope of this paper; we provide here a short summary. The objective of primary (droop) control is to stabilize the network and establish a proportional sharing of load among the DGs.
In islanded operation, inverters are controlled as grid-forming Voltage Source Inverters, having controlled frequencies and voltage magnitudes \cite{QCZ-TH:13}. To accomplish this, a foundation of control loops must be established to regulate the current, voltage, and output impedance of the inverter. This is achieved via an inner control loop for the current, an outer control loop for the voltage, and a virtual impedance loop ensuring the output impedance is desirable at the line frequency (lower half of Figure \ref{Fig:InvTreeNet}) \cite{JMG-LG-JM-MC-JM:05,QCZ-TH:13}.
%
%
%
The reference inputs for the voltage loop are supplied by \emph{droop controllers}, which are heuristic controllers based on active/reactive power decoupling \cite{MCC-DMD-RA:93,JMG-JCV-JM-LGDV-MC:11,JAPL-CLM-AGM:06,MNM-JJW-AK:07,QCZ:13}. For inductive lines, the controllers specify the inverter frequencies $\omega_i$ and voltage magnitudes $E_i$ by\footnote{{Without loss of generality, in islanded mode we consider the droop equations \eqref{eq: droop} without power set points; if desired these can be included as $\omega_i = \omega^* - m_i(P_i-P_{i,\mathrm{set}})$ and $E_i = E^* - n_i(Q_i-Q_{i,\mathrm{set}})$.}}
%
\begin{subequations}
\label{eq: droop}
\begin{align}
	\omega_i &= \omega^* - m_iP_i\,,
	\label{Eq:DroopActive}\\
	\Bigl.
	E_i &= E^* - n_iQ_i\,,
	\label{Eq:DroopReactive}
\end{align}
\end{subequations}
where $\omega^*$ (resp. $E^*$) is a nominal network frequency (resp. voltage), and $P_i$ (resp. $Q_i$) is the measured active (resp. reactive) power injection. The gains $m_i, n_i$ are the \emph{droop coefficients}. 
In \cite{JWSP-FD-FB:12u} a large-signal stability analysis of \eqref{eq: power flow -- active}--\eqref{Eq:DroopActive} was completed, yielding the steady-state network frequency
\begin{equation}\label{Eq:OmegaSS}
\omega_{\rm ss} = \omega^* + \frac{P_0}{\sum_{i=1}^n \frac{1}{m_i}}\,,
\end{equation}
where $P_0$ is the total active power load in the microgrid. Note that the steady-state frequency \eqref{Eq:OmegaSS} is different from the nominal $\omega^*$. Large-signal stability analysis of the voltage droop controller \eqref{Eq:DroopReactive} remains an open problem; see \cite{JWSP-FD-FB:13h} for results on a related droop controller. For non-inductive lines, the appropriate droop controllers take other forms \cite{QCZ-TH:13}.  

%
%
%

\subsection{Review of Secondary Control}
\label{Subsection: Review of Secondary Control}

The removal of the steady-state frequency and voltage deviations generated by the droop controllers \eqref{Eq:DroopActive}--\eqref{Eq:DroopReactive} is accomplished by ``secondary'' integral controllers. 
\smallskip
\subsubsection{Frequency Regulation}
\label{Subsec:FreqReview}

Many techniques have been suggested to restore the network frequency, ranging on the spectrum from centralized to decentralized \cite{JAPL-CLM-AGM:06}, and each with its own advantages and disadvantages. 
{ One centralized technique is to mimic Automatic Generation Control from bulk power systems. This is implemented using area control errors on slow time-scales, a centralized integral controller, and one-to-all communication \cite{JM-JWB-JRB:08}. 
However, this centralized approach conflicts with the microgrid paradigm of distributed generation and autonomous management.}
%
%
A decentralized technique is to use a slower integral control \emph{locally} at each inverter \cite{MCC-DMD-RA:93}. 
This implicitly assumes that the measured local frequency is equal to the steady-state network frequency, and therefore relies on a \emph{separation of time-scales} between the fast, synchronization-enforcing primary droop controller and the slower, secondary integral controller \cite{MCC-DMD-RA:93,JMG-JCV-JM-MC-LGDV:09}. Except in special cases, this approach destroys the power sharing properties established by primary control \cite{MA-DVD-HS-KHJ:14b}, and is too slow to dynamically regulate the grid frequency during rapid load changes.

In \cite{QS-JMG-JCV:14,QS-CS-TD-PP-JCV-JMG:14} control strategies were proposed in which DG units communicate their frequencies, voltages and reactive power injections to one another in order to perform secondary control and share active and reactive power. The methods have two drawbacks: first, all inverters must communicate with all other inverters, requiring a dense communication architecture. 
Second, the controller  gains must be finely tuned in order to maintain active power sharing; see \cite{HB-JWSP-FD-FB:13j} for a detailed analysis.


\subsubsection{Voltage Regulation}
\label{Subsec:VoltReview}
{In high-voltage networks, the sharing among generators of reactive power demand is usually not a major concern due to capacitive compensation of both loads and transmission lines; voltages at generators are therefore regulated to fixed values by the excitation system \cite{JM-JWB-JRB:08}. Voltage regulation has subsequently been adopted as the standard for voltage secondary control in microgrids \cite{MCC-DMD-RA:93, JMG-JCV-JM-LGDV-MC:11}. However, in small-scale microgrid applications, the low ratings of DG units, small electrical distances between units, and the lack of static compensation requires an accurate sharing of reactive power demand among DGs to prevent overloading.}
In Section \ref{Sec:ReactiveSharing} we highlight how voltage regulation and reactive power sharing are conflicting objectives.



Due to the line impedance effect, the voltage droop controller \eqref{Eq:DroopReactive} is unable to share reactive power demand among even identical inverters operating in parallel \cite{JAPL-CLM-AGM:06}. In \cite{QCZ:13}, an alternative primary droop controller was proposed for reactive power sharing among parallel inverters with the same rated voltages. The method requires each unit to have a measurement of the common load voltage, limiting its applicability in complex microgrid scenarios. Similarly, the centralized secondary control architecture proposed in 
\cite{AM-MA-CSS-JMG-JCV:14}
for reactive power sharing requires each unit to communicate with a central controller. The distributed voltage controller proposed in \cite{QS-JMG-JCV:14,QS-CS-TD-PP-JCV-JMG:14}
require all DGs to communicate with all others directly. Moreover, since the controller regulates DG voltages to their nominal values, it is be unable to share reactive power between heterogeneous units connected through varying line impedances. See \cite{AM-MA-CSS-JMG-JCV:14,QS-JMG-JCV:14} and the references therein for more.

\section{Fundamental Limitations of Voltage Control}
\label{Sec:ReactiveSharing}

\begin{figure}[t]
\begin{center}
\includegraphics[width=0.7\columnwidth]{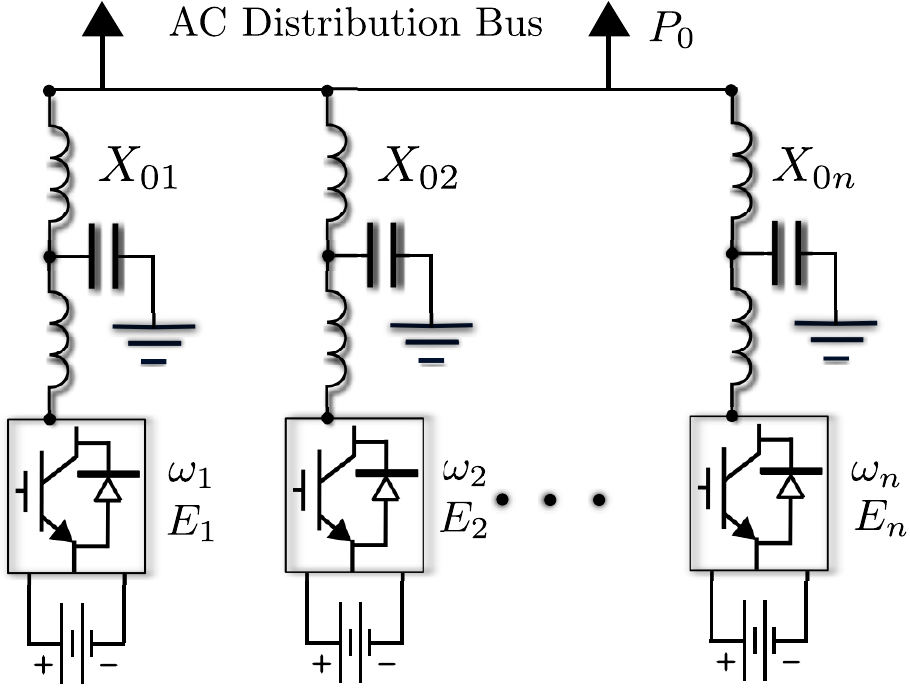}
\captionsetup{justification=raggedright,singlelinecheck=false}
\caption{Schematic of DGs operating in a parallel microgrid.
}
\label{Fig:InvNet}
\end{center}  
\end{figure}

In this section we illustrate the fundamental conflict between two secondary control goals:  voltage regulation and reactive power sharing.
{ For simplicity of exposition, we focus on a parallel microgrid consisting of two identical DGs connected to a common distribution bus (Figure \ref{Fig:InvNet}). The reactances of the two lines connecting the DGs to the common bus are different; in particular, $X_{01} > X_{02}$. }

\begin{figure}[h!]
\begin{center}
\includegraphics[width=0.73\columnwidth]{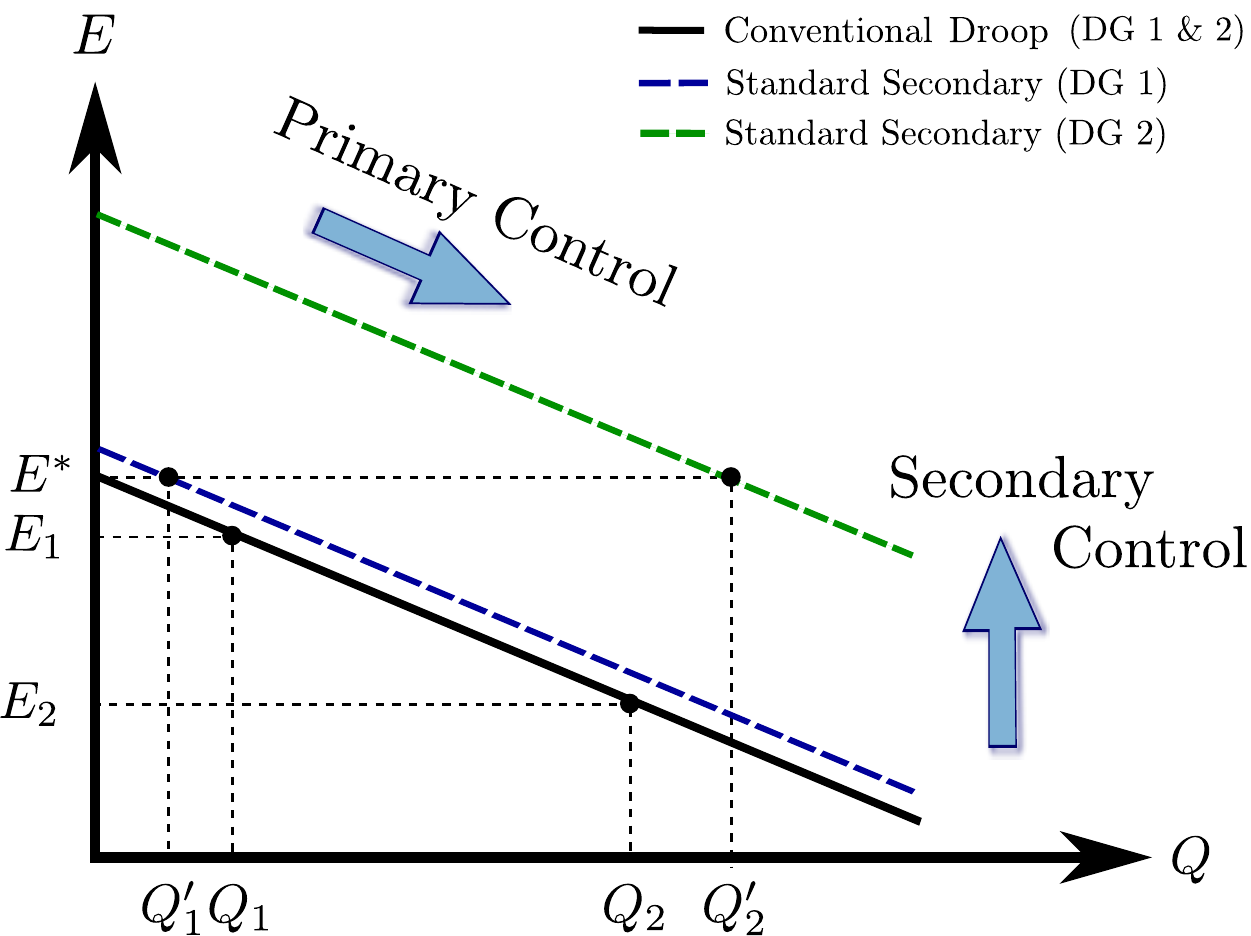}
\captionsetup{justification=raggedright,singlelinecheck=false}
\caption{E-Q droop and standard secondary control for two parallel inverters with identical ratings, operating through reactive lines with $X_{01} > X_{02}$.
}
\label{Fig:VoltChar}
\end{center}  
\end{figure}

Figure \ref{Fig:VoltChar} depicts the E-Q droop characteristics before and after a standard, voltage-regulating secondary control action. Without secondary control, the inverters operate at voltages $E_1$ and $E_2$ with reactive power injections $Q_1$ and $Q_2$ (solid black line). Since $Q_1 \neq Q_2$, reactive power is not shared; this is the ``line impedance effect''. Application of voltage-regulating secondary control ensures that both DG voltage magnitudes are restored to the common rating $E^*$ (dotted blue and green lines are the post-secondary control droop characteristics). Note however that the inverter power injections change to $Q_1^\prime < Q_1$ and $Q_2^\prime > Q_2$. The application of standard secondary control therefore \emph{worsens} the already poor sharing of reactive power between the DGs.

\begin{figure}[t!]
\begin{center}
\includegraphics[width=0.73\columnwidth]{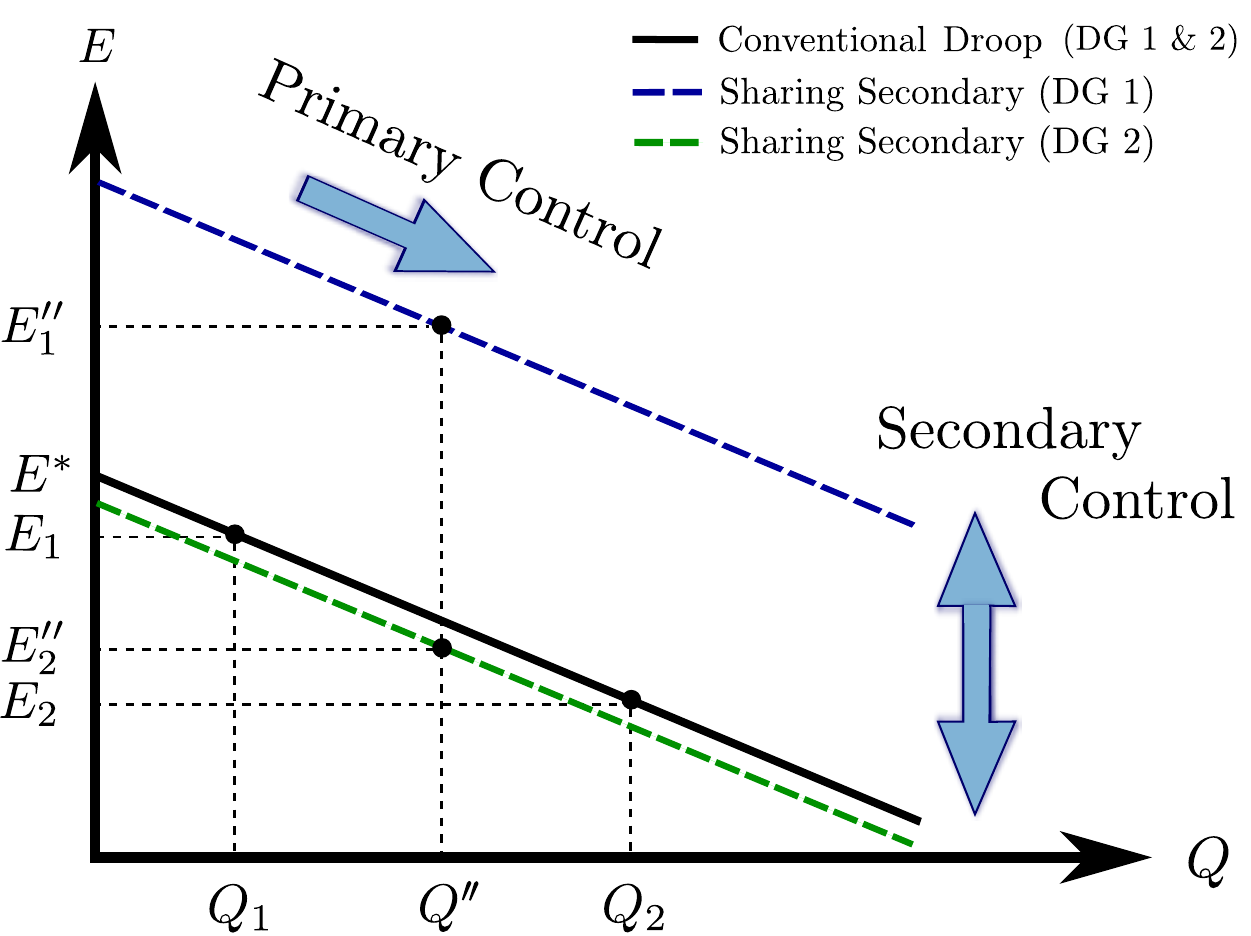}
\captionsetup{justification=raggedright,singlelinecheck=false}
\caption{E-Q droop and power sharing secondary control for two parallel inverters with identical ratings, operating through reactive lines with $X_{01} > X_{02}$.
}
\label{Fig:VoltCharSharing}
\end{center}  
\end{figure}

For the same problem setup, Figure \ref{Fig:VoltCharSharing} depicts the E-Q droop characteristics before and after a power sharing enforcing secondary control action is taken.\footnote{This control action is not uniquely determined; there are many shiftings of the droop characteristics which lead to power sharing (Section \ref{Section:QDAPI}).} While the identical inverters now proportionally share the reactive power by both injecting $Q^{\prime\prime}$, the resulting voltage values $E_1^{\prime\prime}$ and $E_2^{\prime\prime}$ are more dissimilar than they were with only primary control.
\begin{table}[h!]
\begin{center}
\captionsetup{justification=centering,singlelinecheck=false}
\caption{RELATIONSHIPS BETWEEN VOLTAGE MAGNITUDES AND REACTIVE POWER INJECTIONS FOR DIFFERENT CONTROL ACTIONS}
{\renewcommand{\arraystretch}{1}
\begin{tabular}{l|lc}
\toprule
Control Method & Voltage Magnitudes & Reactive Powers\\
\midrule
Primary Control & $E_2 < E_1 < E^*$ & $Q_1 < Q_2$\\
Standard Sec. (${}^\prime$) & $E_2^\prime = E_1^\prime = E^*$ & $Q_1^\prime < Q_1 < Q_2 < Q_2^\prime$ \\
Power Sharing (${}^{\prime\prime}$)& No Relationship & $Q_1^{\prime\prime} =  Q_2^{\prime\prime} = Q^{\prime\prime}$\\
\bottomrule
\end{tabular}
}
\label{Tab:Reactive}
\end{center}
\end{table}

Table \ref{Tab:Reactive} collects the relationships between voltage magnitudes and reactive power injections for the different control actions described above. 
We observe that \textemdash{} except under special circumstances \textemdash{} precise voltage regulation leads to large errors in reactive power sharing, as shown in Figure \ref{Fig:VoltChar}. Conversely, the objective of reactive power sharing does not uniquely determine the DG bus voltages, and when implemented improperly can result in poor voltage profiles as shown in Figure \ref{Fig:VoltCharSharing}. The accuracy of reactive power sharing that can be achieved therefore depends on both the upper and lower limits for the DG voltage magnitudes, and on the homogeneity of the line reactances.
%
%
We conclude that an ideal secondary voltage controller should allow for a \emph{tunable compromise} between voltage regulation and reactive power sharing. 

\section{Distributed Averaging Proportional Integral (DAPI) Controllers for Microgrids}
\label{Section:DAPI}

As mentioned in Section \ref{Section: Introduction}, communication has been identified as an essential ingredient for high-performance secondary control. We now introduce the DAPI controllers, which combine droop and integral control with averaging algorithms from multi-agent systems \cite{WR-RWB-EMA:07}. To build intuition for our designs, we first briefly review continuous-time averaging over networks.

\subsection{Review of Continuous-Time Distributed Averaging}
\label{Section:ReviewofAveraging}

\begin{figure}[h!]
\begin{center}
\includegraphics[width=0.83\columnwidth]{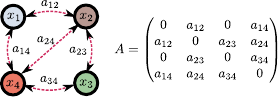}
\captionsetup{justification=raggedright,singlelinecheck=false}
\caption{Example of adjacency matrix construction for four DGs, with $\mathcal{V} = \{1,2,3,4\}$ and $\mathcal{E} = \{(1,2),(2,1),(2,3),(3,2),(3,4),\ldots\}$.}
\label{Fig:Averaging}
\end{center}  
\end{figure}

The communication layer between DG's will be described by a weighted graph $G(\mathcal{V},\mathcal{E},A)$ where $\mathcal{V} = \{1,\ldots,n\}$ is a labeling of the DGs, $\mathcal{E} \subseteq \mathcal{V} \times \mathcal{V}$ is the set of communication links, and $A$ is the $n \times n$ \emph{weighted adjacency matrix} of the graph, with elements $a_{ij} = a_{ji} \geq 0$. In particular, one writes that $(i,j) \in \mathcal{E}$ if node $i$ sends information directly to node $j$, and in this case $a_{ji} > 0$. Thus, the sparsity pattern of the adjacency matrix $A$ encodes the topology of the communication layer (Figure \ref{Fig:Averaging}).
If to each node $i \in \{1,\ldots,n\}$ we assign a scalar value $x_i$, a commonly studied update rule is for node $i$ to adjust its value $x_i$ according to
\begin{equation}\label{Eq:Consensus}
\dot{x}_i = -\sum_{j=1}^n \nolimits a_{ij}(x_i-x_j)\,.
\end{equation}
Equation \eqref{Eq:Consensus} is called continuous-time distributed averaging, or ``consensus''. To interpret \eqref{Eq:Consensus}, define the convex weights\footnote{The weights are convex because $w_{ij} \geq 0$ and $\sum_{j=1}^n w_{ij} = 1$.} $w_{ij} = a_{ij}/(\sum_{k=1}^n a_{ik})$, and rearrange \eqref{Eq:Consensus} to obtain
\begin{equation}\label{Eq:RearrangeConsensus}
\frac{1}{\sum_{j=1}^n \nolimits a_{ij}}\dot{x}_i = -x_i + \sum_{j=1}^n \nolimits w_{ij}x_j\,.
\end{equation}
Thus, with time-constant $1/\sum_{j=1}^n a_{ij}$, the variable $x_i$ evolves toward a weighted average of its neighbors values $x_j$, with averaging weights $w_{ij}$. If the communication network $G(\mathcal{V},\mathcal{E},A)$ is connected, this dynamic process results in all variables $x_i$ \emph{converging to a common value} $x_i = x_j = $ const. (see Remark \ref{Remark:CommunicationRequirements}) \cite{ROS-JAF-RMM:07, WR-RWB-EMA:07}. We now apply these ideas from continuous-time distributed averaging to microgrid control.
%
%

\subsection{Frequency Regulation and Active Power Sharing}
\label{Section:PDAPI}

We propose the \emph{distributed-averaging proportional-integral (DAPI) frequency controller}
\begin{subequations}\label{Eq:DAPIActive}
\begin{align}
\omega_i &= \omega^* - m_iP_i + \Omega_i\,,
\label{Eq:Primary}
\\
k_i\frac{\mathrm{d}\Omega_i}{\mathrm{d}t} &= -(\omega_i-\omega^*) - \sum_{j=1}^n a_{ij}\left(\Omega_i-\Omega_j\right)\,,
\label{Eq:Secondary}
\end{align}
\end{subequations}
where $\Omega_i$ is the secondary control variable and $k_{i}$ is a positive gain. The first equation \eqref{Eq:Primary} is the standard droop controller with the additional secondary control input $\Omega_i$. To understand the second equation, we consider two cases.

\textbf{Case 1 ($\boldsymbol{A = 0}$): } In this case, there is no communication among DGs, and $\Omega_i$ is the integral of the local frequency error $\omega_i-\omega^*$ with gain $1/k_i$. In steady-state, the derivative on the left of \eqref{Eq:Secondary} is zero, and hence $\omega_i = \omega^*$ for each DG $i$. That is, the network frequency has been regulated. However, depending on the initial conditions and controller gains, the variables $\Omega_i$ may converge to different values, and shift their respective droop curves by different amounts. This unwanted degree of freedom leads to poor active power sharing.

\begin{figure}[t]
\begin{center}
\includegraphics[width=0.75\columnwidth]{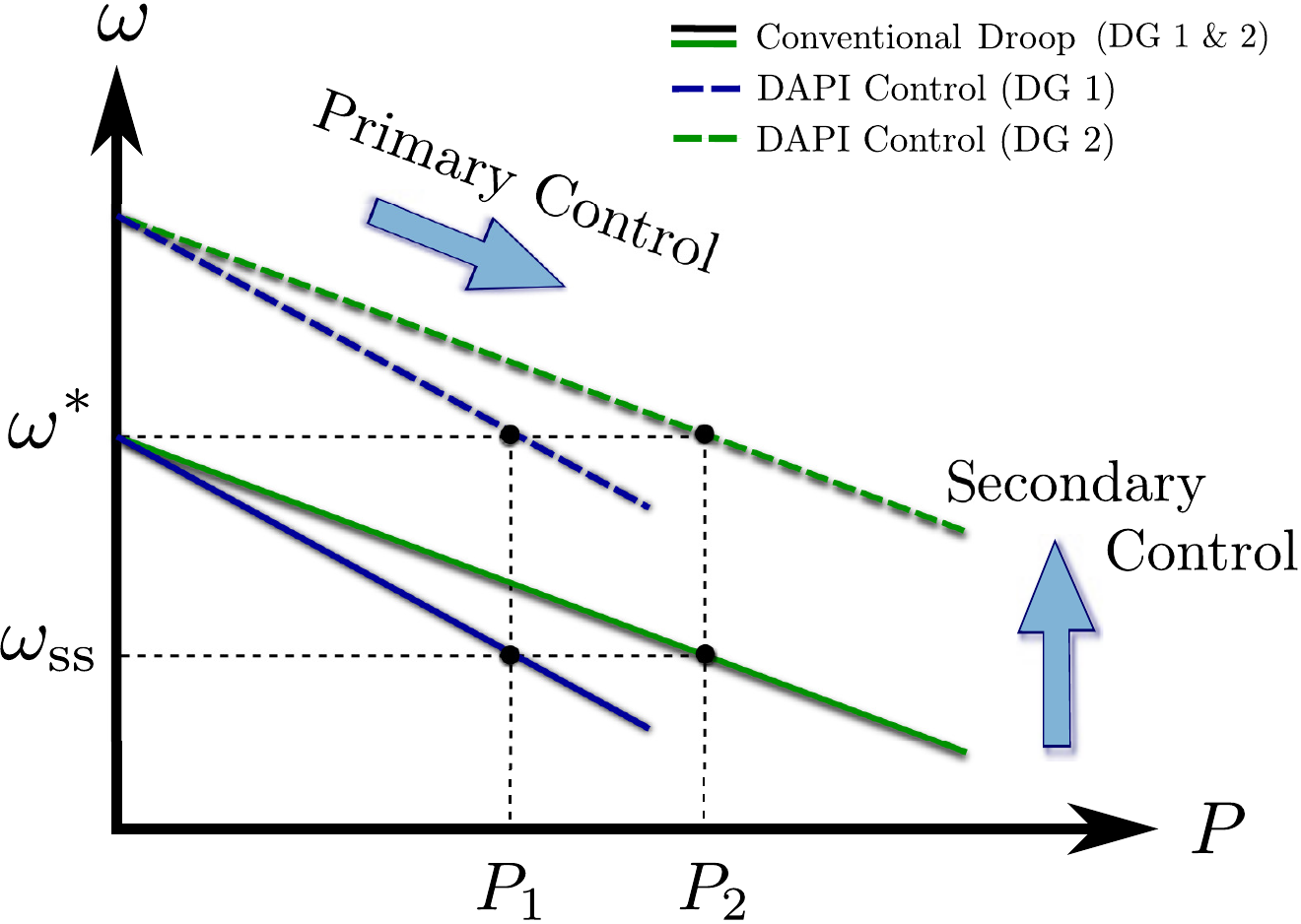}
\captionsetup{justification=raggedright,singlelinecheck=false}
\caption{Droop characteristics with (dashed lines) and without (solid lines) secondary control action. One can interpret the DAPI secondary control action as a uniform shifting of all droop characteristics by an amount $\omega^*-\omega_{\rm ss}$.
}
\label{Fig:FreqChar}
\end{center}  
\end{figure}

\textbf{Case 2 ($\boldsymbol{A \neq 0}$): } Now, consider the case where we include the diffusive averaging terms $a_{ij}(\Omega_i-\Omega_j)$. As before, in steady-state the derivative on the left-hand side of \eqref{Eq:Secondary} must be zero, and hence $\omega_i=\omega^*$. However, from the discussion in Section \ref{Section:ReviewofAveraging}, we also must have $\Omega_i = \Omega_j$ for all inverters $i,j$. That is, the DGs must agree on how much to shift the droop characteristics.
This ensures that all droop curves are shifted by the same amount equal to $\omega^*-\omega_{\rm ss}$ (Figure \ref{Fig:FreqChar}), guaranteeing active power sharing is maintained. This performance is not dependent on the controller gains $k_i$ and $a_{ij}$, which determine only the transient behavior of the controller (see Table \ref{Tab:Tune}).

\begin{remark}\label{Remark:CommunicationRequirements}\textbf{(Communication Requirements for DAPI Control).} DAPI control requires that neighboring DG units exchange information to collectively perform secondary control. To ensure power sharing among all units, the communication network among DGs must be \emph{connected}: there must be a path in the communication graph between any two nodes, as in Figure \ref{Fig:Averaging}. While here we consider the controller in continuous-time with bi-directional communication, our assumptions can be relaxed to allow for asymmetric, asynchronous and discrete-time communication with delays \cite{WR-RWB-EMA:07,ROS-JAF-RMM:07}.
\end{remark}


\subsection{Voltage Regulation and Reactive Power Sharing}
\label{Section:QDAPI}

As noted in Section \ref{Subsection: Review of Secondary Control}, the E-Q droop controller \eqref{Eq:DroopReactive} is unable to share reactive power between DGs. Moreover, in Section \ref{Sec:ReactiveSharing} we described the conflict between voltage regulation and reactive power sharing.
With these problems in mind, we propose the second DAPI controller
\begin{subequations}\label{Eq:DAPIReactive}
\begin{align}
E_i &= E^* - n_iQ_i + e_i\,,
\label{Eq:PrimaryReactive}
\\
\kappa_i\frac{\mathrm{d}e_i}{\mathrm{d}t} &= -\beta_i(E_i-E^*)-\sum_{j=1}^n b_{ij}\left(\frac{Q_i}{Q_i^*}-\frac{Q_j}{Q_j^*}\right)\,,
\label{Eq:SecondaryReactive}
\end{align}
\end{subequations}
where $e_i$ is the secondary control variable, $Q_i^*$ is the $i$th DGs reactive power rating, and $\beta_i, \kappa_i$ are positive gains. The $n \times n$ matrix $B$ with elements $b_{ij} > 0$ is the adjacency matrix of a  communication network between the DGs. The secondary controller \eqref{Eq:SecondaryReactive} achieves a tunable compromise between voltage regulation and reactive power sharing. We consider four cases:

\textbf{Case 1 ($\boldsymbol{\beta = 0}$, $\boldsymbol{B \neq 0}$): }
In this case the first term in \eqref{Eq:SecondaryReactive} is disappears, leaving only the second term. Steady-state requires the derivative on the left-hand side of \eqref{Eq:SecondaryReactive} to be zero, which occurs if and only if $Q_i/Q_i^* = Q_j/Q_{j}^*$ for all inverters. Thus, the steady-state is a power sharing configuration. The secondary control variables $e_i$ converge to values which shift the individual droop curves as necessary to establish proportional power sharing, see Figure \ref{Fig:VoltCharSharing}. However, as discussed in Section \ref{Sec:ReactiveSharing}, under such a control action DG voltages can drift quite far from their nominal values. 


\textbf{Case 2 ($\boldsymbol{\beta \neq 0}$, $\boldsymbol{B = 0}$):}
In this case the second term in \eqref{Eq:SecondaryReactive} disappears, and the controller reduces to the standard decentralized voltage-regulating secondary control discussed in Section \ref{Sec:ReactiveSharing}. Reactive power is shared poorly (Figure \ref{Fig:VoltChar}).

\textbf{Case 3 ($\boldsymbol{\beta \neq 0}$, $\boldsymbol{B \neq 0}$):} In this regime \eqref{Eq:PrimaryReactive}--\eqref{Eq:SecondaryReactive} achieves a compromise between reactive power sharing and voltage regulation based on the relative sizes of the gains $\beta_i$ and $b_{ij}$.

\textbf{Case 4 (Smart Tuning): } As a specialization of Case 3, consider having a specific DG $i$ implement the controller \eqref{Eq:SecondaryReactive} with $\beta_i \neq 0$ and $b_{ij} = 0$, while the all other DGs $j\neq i$ implement \eqref{Eq:SecondaryReactive} with $\beta_j = 0$ and $b_{jk} \neq 0$.\footnote{This directed communication tuning requires that DG $i$ sends information to at least one neighbor.} That is, DG $i$ regulates its voltage to the nominal value, and the voltages at DGs $j\neq i$ are then controlled to share power in a manner consistent with the voltage regulation of DG $i$ (cf. Section \ref{Sec:ReactiveSharing}). This tuning sets up a ``leader-follower'' \cite{FB-JC-SM:09} relationship among the DGs, where the voltages at DGs $j \neq i$ will form a cluster around the voltage value of $E_i = E^*$ of DG $i$.

\smallskip

The above cases are tested experimentally in Section \ref{Sec:Study1}.

\smallskip

\begin{remark}\label{Rem:Interp}\textbf{(Remarks on DAPI Control).}
The communication layers between DG units described the adjacency matrices $A$ and $B$ are \emph{design variables of the DAPI controllers}. This customizable architecture allows for design flexibility. For example, to add redundancy against communication channels being permanently disconnected, supplementary communication can be introduced. Note that the communication architecture need not mirror the electrical topology of the network (Figure \ref{Fig:Schematic}), and that the controllers do not rely on high-gain techniques such as feedback linearization \cite{AB-AD-FLL-JMG:13}. {A detailed schematic of the DAPI control architecture is shown in Figure \ref{Fig:InvTreeNet}.}

\begin{figure}[t]
\begin{center}
\includegraphics[width=0.9\columnwidth]{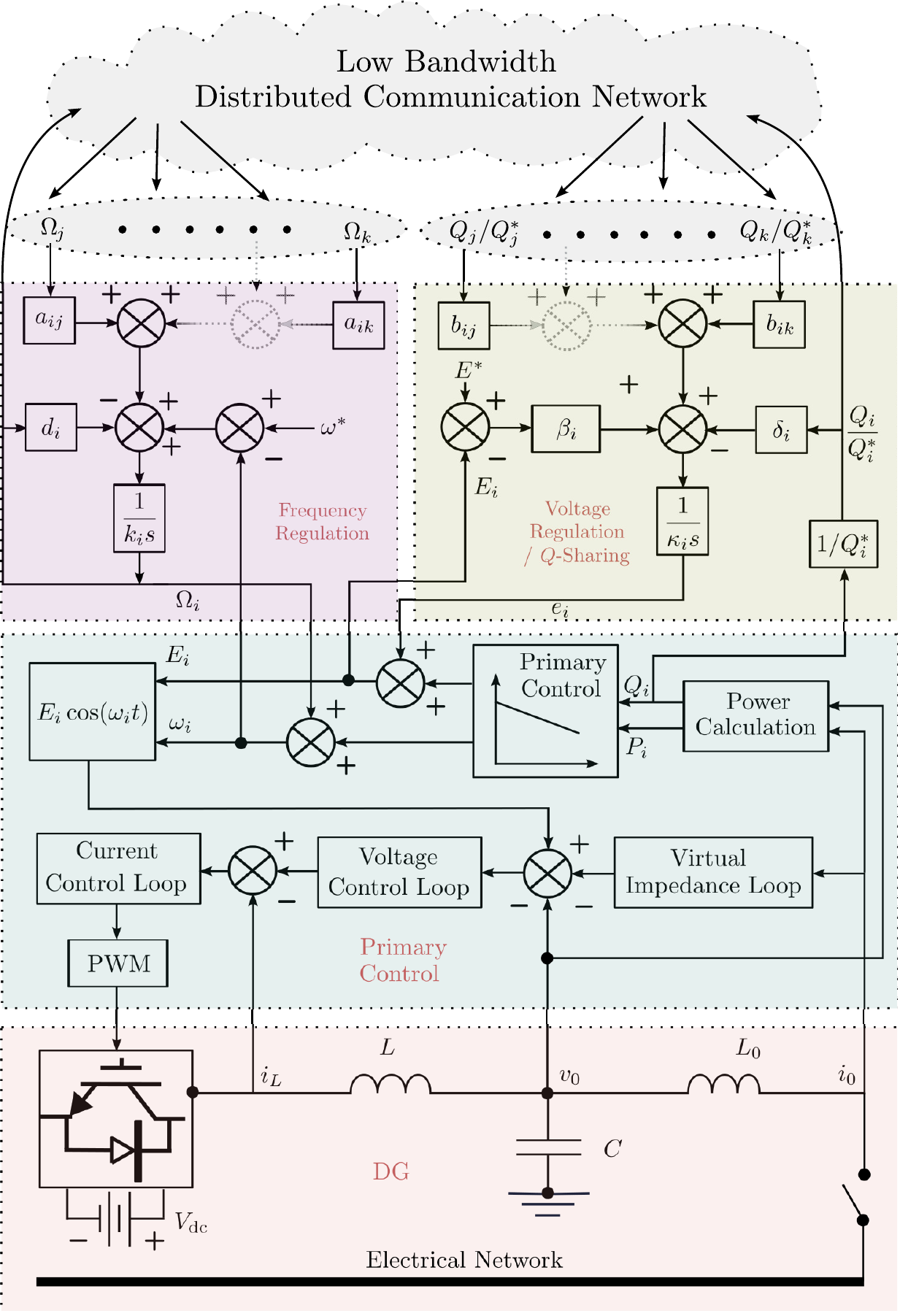}
\captionsetup{justification=raggedright,singlelinecheck=false}
\caption{Block diagram of proposed control architecture for a single DG. For simplicity we have have abbreviated $d_i = \sum_{j=1}^n a_{ij}$ and $\delta_i = \sum_{j=1}^n b_{ij}$.}
\label{Fig:InvTreeNet}
\end{center}  
\end{figure}


{ The time-constants $k_i$ and $\kappa_i$ in \eqref{Eq:Secondary} and \eqref{Eq:SecondaryReactive} allow for a precise tuning of the secondary control speed. A conventional choice is to make $k_i$ and $\kappa_i$ sufficiently large, enforcing a time-scale separation between primary and secondary control. This however is not \emph{required} -- our experimental results suggest that primary and secondary control can be performed on similar time scales without stability issues or performance degradation.\footnote{In fact, in \cite{JWSP-FD-FB:12u} it was shown that the frequency controller \eqref{Eq:Primary}--\eqref{Eq:Secondary} is stabilizing for any choice of gains $k_i$. }  Table \ref{Tab:Tune} provides a simple qualitative reference for the effects of the control parameters in \eqref{Eq:DAPIActive} and \eqref{Eq:DAPIReactive}.}

\begin{table}[t]
\begin{center}
\captionsetup{justification=centering,singlelinecheck=false}
\caption{{QUALITATIVE EFFECTS OF CONTROLLER GAINS}}
{\renewcommand{\arraystretch}{1}
\begin{tabular}{l|l}
Gain & Qualitative Change Upon Increase\\
\midrule
$k_i$ & Slows frequency regulation at DG $i$\\
$\kappa_i$ & Slows voltage regulation / $Q$-sharing at DG $i$\\
$a_{ij}$ & Faster $P$-sharing between DGs $i$ and $j$\\
$b_{ij}$ & Improved steady-state $Q$-sharing between DGs $i$ and $j$\\
$\beta_i$ & Improved steady-state voltage regulation at DG $i$\\
\bottomrule
\end{tabular}
}
\label{Tab:Tune}
\end{center}
\end{table}
%

%
%

\end{remark}

{ 
\section{Stability \& Performance of DAPI Control}
\label{Sec:Stability}

A large-signal nonlinear stability analysis of the frequency DAPI controller \eqref{Eq:DAPIActive} can be found in \cite{JWSP-FD-FB:12u}. 
While the secondary frequency controller \eqref{Eq:Secondary} will never destabilize the primary controller \eqref{Eq:Primary}, the secondary voltage controller \eqref{Eq:SecondaryReactive} can potentially destabilize \eqref{Eq:PrimaryReactive}. 
This possibility exists due to the previously discussed conflict between reactive power sharing and voltage regulation. 
A full nonlinear stability analysis of the voltage/reactive power DAPI controller \eqref{Eq:PrimaryReactive}--\eqref{Eq:SecondaryReactive} is extremely challenging and beyond the scope of this article; a partial analysis for a simpler controller can be found in \cite{JC-TS-JR-TS:14}. 
{In Section \ref{Sec:SS} we present a small-signal stability analysis of \eqref{Eq:PrimaryReactive}--\eqref{Eq:SecondaryReactive}, along with sufficient conditions which ensure stable operation. In Section \ref{Sec:Trans} we explore the effect of the controller gains in \eqref{Eq:DAPIActive},\eqref{Eq:DAPIReactive} on the transient performance of the closed-loop system.}

\subsection{Small-Signal Stability of Voltage DAPI Control}
\label{Sec:SS}

To avoid unnecessary technical complications, we model any delay in adjusting the output voltage in \eqref{Eq:PrimaryReactive} with a simple low-pass filter, yielding the dynamic system $\mathrm{d}E_i/\mathrm{d}t = -(E_i-E_i^*) - n_iQ_i + e_i$, and assume loads are impedances collocated with DGs. 
Both of these assumptions can be relaxed at the expense of more complicated formulae. 
Under the standard decoupling assumption in which reactive power is related strongly to differences in voltage magnitudes \cite{JM-JWB-JRB:08}, the reactive power injection \eqref{eq: power flow -- reactive} at the $i$th DG takes the form
\begin{equation}\label{Eq:LoadFlow}
Q_i = -E_i^2Y_{\mathrm{load},ii} + E_i\sum_{j=1}^{n} \nolimits Y_{\mathrm{bus},ij}(E_i-E_j)\,,
\end{equation}
%
%
%
%
where $Y_{\mathrm{load}}$ is diagonal matrix of load susceptances and $Y_{\rm bus} = Y_{\rm bus}^T$ is the microgrid's bus admittance matrix \cite{JM-JWB-JRB:08}. In vector notation, the system equations \eqref{Eq:DAPIReactive},\eqref{Eq:LoadFlow} are\footnote{Here $[z]$ denotes the diagonal matrix with the vector $z$ along the diagonal.}%
\begin{subequations}\label{Eq:DAPIReactiveV}
\begin{align}
\dot{E}  &= -(E - E^*) - NQ + e\,,
\label{Eq:PrimaryReactiveV}
\\
\kappa\dot{e} &= -\beta(E-E^*)-L_c[Q^*]^{-1}Q\,,
\label{Eq:SecondaryReactiveV}
\\
Q &= 
[E]YE\,,
\label{Eq:LoadFlowV}
\end{align}
\end{subequations}
where $E, E^*$ and $e$ are the vectors of voltage magnitudes, voltage set points, and secondary control variables, $N, \beta$ and $\kappa$ are diagonal matrices of controller gains, $Q$ and $Q^*$ are the vectors of DG reactive power injections and reactive power ratings, $L_{\rm c} = \mathrm{diag}(\sum_{j=1}^n b_{ij})-B$ is the Laplacian matrix \cite{ROS-JAF-RMM:07} corresponding to the communication network among the DGs, and $Y = -(Y_{\rm bus} + Y_{\rm load})$. 
%
%
%
%
When implementing the controller \eqref{Eq:DAPIReactive} in practice, the voltages $E_i$ will remain near their nominal values $E^*$. We can exploit this to obtain a linear dynamic system by making the approximation that $[E] \simeq [E^*]$ in \eqref{Eq:LoadFlowV}; details on this approximation technique can be found in \cite{BG-JWSP-FD-SZ-FB:13zb}.
After making this approximation and inserting \eqref{Eq:LoadFlowV} into \eqref{Eq:SecondaryReactiveV}, the nonlinear system \eqref{Eq:PrimaryReactiveV}--\eqref{Eq:SecondaryReactiveV} becomes the linear system 
\begin{equation}\label{Eq:LinSys}
\dot{x} = Wx + u\,,
\end{equation}
where $x = (E,e)$, $u = (E^*,\kappa^{-1}\beta E^*)$, and
%
%
\begin{equation*}
W = \begin{pmatrix}
-W_{1} & I_n \\
-W_{2} & 0_n
\end{pmatrix} = \begin{pmatrix}
-(I_n + N[E^*]Y) & I_n\\
-\kappa^{-1}(\beta+L_c[Q^*]^{-1}[E^*]Y) & 0_n
\end{pmatrix}\,,
\end{equation*}
where $I_n$ (resp. $0_n$) is the $n \times n$ identity matrix (resp. zero matrix). 
For future use, we note that all eigenvalues of $-W_1$ are real and negative since $W_1$ is similar to a symmetric $M$-matrix, as can be verified by using the similarity transform $TW_1T^{-1}$ where $T = N_I^{-1/2}[E_I^*]^{-1/2}$.
%
%
%
%
%
%
%
%
%
%
%
%
%
%
%
%
%
{ We now derive sufficient conditions under which the linearized system \eqref{Eq:LinSys} is exponentially stable. Specifically, we will assume that the characteristic polynomial $\det(sI_{2n}-W) = 0$ of \eqref{Eq:LinSys} has a root in the closed right-half complex plane, and derive conditions under which this assumption is contradicted. These conditions will therefore ensure all characteristic roots are in the left-half complex plane, and thus ensure stability.}
%
%
Since $-W_1$ has negative eigenvalues, it follows that $\det(s I_n + W_1) \neq 0$, and using the Schur complement determinant formulae for block matrices we may simplify the characteristic polynomial as
%
\begin{equation}\label{Eq:CharPoly}
\det(sI_{2n} - W) = \det(s^2 I_n + s W_1 + W_2) = 0\,.
\end{equation}
Since the determinant is zero, the matrix $\,s^2 I_n + s W_1 + W_2$ must be singular, and therefore the polynomial \eqref{Eq:CharPoly} has a solution if and only if $x^T(s^2I_n+s W_1+W_2)x = 0$ for some real vector $x$ of unit length. The latter is simply a scalar quadratic equation of the form $\,s^2 + \alpha_1 s+ \alpha_2 = 0$, where $\alpha_1 = x^TW_1x$ and $\alpha_2 = x^TW_2x$. If it is true that
\begin{subequations}\label{Eq:Cond}
\begin{align}\label{Eq:Cond1}
\lambda_{\rm min}(W_1 + W_1^T) &> 0\,,\\
\label{Eq:Cond2}
\lambda_{\rm min}(W_2 + W_2^T) &> 0\,,
\end{align}
\end{subequations}
where $\lambda_{\rm min}(\cdot)$ is the smallest eigenvalue of the matrix argument, then $\alpha_1,\alpha_2 > 0$ and all solutions of $s^2 + \alpha_1s+\alpha_2 = 0$ satisfy $\mathrm{Re}(s) < 0$ by the Routh-Hurwitz criterion.
%
%
%
%
This contradicts our assumption that the characteristic polynomial has a closed right-half plane root, and hence under the conditions \eqref{Eq:Cond} the linearized system is exponentially stable.
%
%
%
%

Let us now physically interpret the stability conditions \eqref{Eq:Cond}. The first condition \eqref{Eq:Cond1} restricts the DGs from being too dissimilar. For example, if all DGs have the same droop gains $n_i$ and voltage set points $E^*$, then $W_1$ is a scalar values times $Y$ and \eqref{Eq:Cond1} is always true. For dissimilar DGs, the intuition for \eqref{Eq:Cond1} is that given equal voltage set points, DGs with high ratings should be connected to the microgrid through stiff lines of high admittance. 
%
%
%
%
%
{ To understand the second condition \eqref{Eq:Cond2}, we first consider the case of pure voltage regulation (Case 2 in Section \ref{Section:QDAPI}) where $\beta_i \neq 0$ and $L_{\rm c} = 0_n$. Then $W_2 = \kappa^{-1}\beta > 0$ is diagonal and \eqref{Eq:Cond2} is satisfied. The voltage regulation gains $\beta_i$ always act to stabilize the system. Since eigenvalues are continuous functions of matrix parameters, the system is also stable for non-zero but sufficiently small power sharing gains $b_{ij}$. In the more general case where the power sharing gains $b_{ij}$ are also non-zero, the condition \eqref{Eq:Cond2} properly accounts for the complicated interplay between the microgrid's electrical stiffness matrix $Y$ and the averaging control action $L_{\rm c}$ in the product $L_c[Q^*]^{-1}[E^*]Y$. Intuitively, \eqref{Eq:Cond2} will be satisfied when all line impedances are sufficiently uniform and all DGs are sufficiently similar, since in this case reactive power sharing is not in strong conflict with the line impedance effect (cf. Section \ref{Sec:ReactiveSharing}). The stability conditions \eqref{Eq:Cond} are both satisfied in all experiments presented in Section \ref{Sec:Experiments}.
%
%
}

%
%
%
%
%
%
%
} 

\subsection{Transient Performance of DAPI Control}
\label{Sec:Trans}

\begin{figure}[t]
\begin{center}
\includegraphics[width=1\columnwidth]{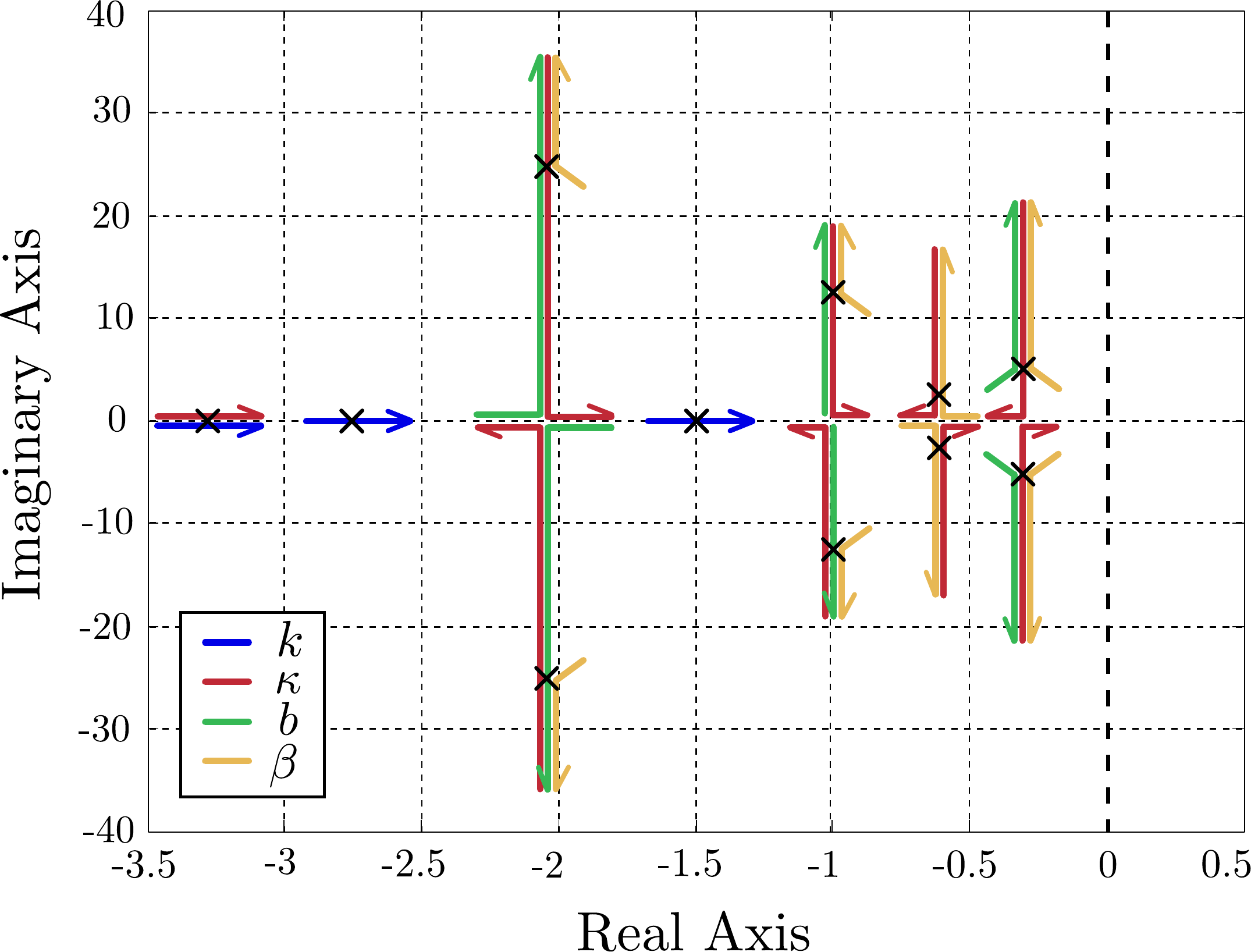}
\captionsetup{justification=raggedright,singlelinecheck=false}
\caption{{ Eigenvalue traces of closed-loop system \eqref{eq: power flow},\eqref{Eq:DAPIActive},\eqref{Eq:DAPIReactive} as controller gains are varied. Arrows indicate the direction increasing gain. The absence of a trace indicates that the parameter under consideration has negligible effect on the respective eigenvalue. System parameters are taken from Study 1c of Section \ref{Sec:Experiments}. Black crosses indicate eigenvalue locations for the nominal gains used in Study 1c. Several fast eigenvalues are omitted for clarity.}}
\label{Fig:rlocus}
\end{center}  
\end{figure}

{ The impact of the controller gains $k_i, a_{ij}, \kappa_i, \beta_i$ and $b_{ij}$ on the steady-state equilibrium was discussed in detail in Section \ref{Section:DAPI} and summarized in Table \ref{Tab:Tune}. 
We now examine the impact of these controller gains on the system's \emph{transient performance}.
To do this, we consider a case study with four DGs (Figure \ref{Fig:Schematic}), with the system parameters of Table \ref{Tab:Parameters}. 
The communication network among the DGs is a ring, and the controller gains $a_{ij}$ and $b_{ij}$ are given by \eqref{Eq:CommTop}; $a_{ij}$ are constants, while $b_{ij}$ are parameterized by a single constant $b$. 
For all DGs $i = 1,\ldots,4$, the other control parameters $k_i, \kappa_i$ and $\beta_i$ are taken as uniform constants $k, \kappa$ and $\beta$, respectively. 
The nominal values for these gains are $k = 1.7\,$s, $\kappa = 1\,$s, $\beta = 1.2$, and $b = 180\,$V (the same as in Study 1c of Section \ref{Sec:Experiments}).

We increment the gains independently in intervals around their nominal values, and for each iteration we (i) numerically determine the system operating point from \eqref{eq: power flow},\eqref{Eq:DAPIActive},\eqref{Eq:DAPIReactive}, (ii) linearize the closed-loop system around the operating point, and (iii) plot the eigenvalues of the linearization. These eigenvalue traces are displayed in Figure \ref{Fig:rlocus}, where black crosses denote the eigenvalues for the nominal controller gains and arrows indicate the direction of increasing gain.

Eigenvalues on the real axis are strongly associated with the frequency dynamics \eqref{Eq:Primary}--\eqref{Eq:Secondary}, while complex conjugate eigenvalues are associated with the voltage dynamics \eqref{Eq:PrimaryReactive}--\eqref{Eq:SecondaryReactive}. 
These conjugate eigenvalues lead to an underdamped voltage response: physically, this is a manifestation of the line impedance effect, which the controller \eqref{Eq:PrimaryReactive}--\eqref{Eq:SecondaryReactive} must overcome to establish reactive power sharing.
As the frequency time-constant $k$ is increased (blue), real eigenvalues move towards the origin leading to slow, smooth frequency/active power response. 
Conversely, decreasing $k$ leads to fast (but still overdamped) frequency regulation. Increasing $a_{ij}$ has an effect nearly identical to decreasing $k$, and we have therefore omitted the trace and held $a_{ij}$ constant. 
Increasing the voltage time-constant $\kappa$ (red) causes the underdamped conjugate eigenvalues to collapse onto the real axis, leading to an overdamped voltage/reactive power response for sufficiently slow secondary control. 
Increasing either feedback gain $b$ or $\beta$ (green and gold) results in an increasingly  underdamped voltage/reactive power response.

Taken together, Table \ref{Tab:Tune}, the stability conditions \eqref{Eq:Cond}, and the eigenvalue traces of Figure \ref{Fig:rlocus} provide a solid foundation for understanding the DAPI controllers \eqref{Eq:DAPIActive}--\eqref{Eq:DAPIReactive}. Our experimental results demonstrate that despite the simplifying assumptions in the preceding analysis, the DAPI controllers \eqref{Eq:DAPIActive}--\eqref{Eq:DAPIReactive} can be tuned for both stability and high performance.}

\section{Experimental Results}
\label{Sec:Experiments}
\begin{figure}[t]
\begin{center}
\includegraphics[width=0.85\columnwidth]{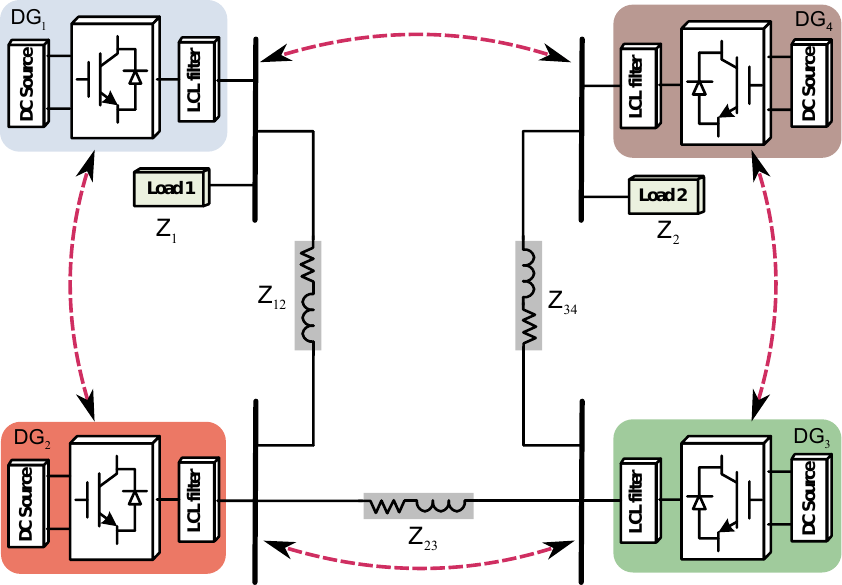}
\captionsetup{justification=raggedright,singlelinecheck=false}
\caption{Schematic of the experimental microgrid setup, consisting of four DGs interconnected through heterogeneous impedances. Loads are collocated at DGs one and four. Red dotted lines denote communication links.}
\label{Fig:Schematic}
\end{center}  
\end{figure}

Experiments were performed at the Intelligent Microgrid Laboratory (Aalborg University, Denmark)  
to validate the DAPI controllers presented in Section \ref{Section:DAPI}. { A schematic of the experimental setup is shown in Figure \ref{Fig:Schematic}, consisting of four DGs interconnected through impedances}. Loads are present locally at units 1 and 4, and units 1 and 4 are rated for twice as much power as units 1 and 3 (Table \ref{Tab:Parameters}). The DAPI controllers \eqref{Eq:Primary}--\eqref{Eq:SecondaryReactive} were implemented in Simulink\textsuperscript{\textregistered}, with measurements recorded via a dSPACE\textsuperscript{\textregistered} 1006. See \cite{JMG-JCV-JM-LGDV-MC:11,QCZ-TH:13} for details on the inner voltage, current and impedance control loops.

\smallskip

This section is organized into four studies, beginning with a characterization of controller performance, and then examining robustness under communication link failure, heterogeneous controller gains, and plug-and-play operation. The communication structure is shown in Figure \ref{Fig:Schematic}, with the adjacency matrices $A = [a_{ij}]$ in \eqref{Eq:Secondary} and $B = [b_{ij}]$ in \eqref{Eq:SecondaryReactive} being
\begin{equation}\label{Eq:CommTop}
A = \begin{pmatrix}
0 & 1 & 0 & 1\\
1 & 0 & 1 & 0\\
0 & 1 & 0 & 1\\
1 & 0 & 1 & 0
\end{pmatrix}\,,\quad B = b\cdot \begin{pmatrix}
0 & 1 & 0 & 1\\
1 & 0 & 1 & 0\\
0 & 1 & 0 & 1\\
1 & 0 & 1 & 0
\end{pmatrix}\,,
\end{equation}
where $b \geq 0$ varies depending on the study under consideration. All other parameters are as reported in Table \ref{Tab:Parameters}. All plots are color-coded in correspondence with Figure \ref{Fig:Schematic}: DG 1 (blue), DG 2 (red), DG 3 (green), and DG 4 (brown).


\begin{table}[t]
\begin{center}
\captionsetup{justification=centering,singlelinecheck=false}
\caption{ELECTRICAL AND CONTROL PARAMETERS}
{\renewcommand{\arraystretch}{1}
\begin{tabular}{lll}
Parameter & Symbol & Value\\
\toprule
\multicolumn{3}{c}{Electrical Setup} \\
\midrule
Nominal Frequency & $\omega^*/2\pi$ & 50\,Hz\\
DC Voltage & $\,V_{\rm dc}$ & 650\,V\\
Nominal Voltages & $E^*$ & {325.3\,V (230\,V rms)}\\
Filter Capacitance & $C$ & 25\,$\mu$F\\
Filter Inductance & $L_f$ & 1.8\,mH\\
Output Impedance & $L_0$ & 1.8\,mH\\
Line Impedance (1,2) & $Z_{12}$ & $R = 0.8\,\Omega$, $L = 3.6\,$mH\\
Line Impedance (2,3) & $Z_{23}$ & $R = 0.4\,\Omega$, $L = 1.8\,$mH\\
Line Impedance (3,4) & $Z_{34}$ & $R = 0.7\,\Omega$, $L = 1.9\,$mH\\
\midrule
\multicolumn{3}{c}{Control Parameters} \\
\midrule
\end{tabular}
\begin{tabular}{llll}
Parameter & Symbol & DGs 1\&4 & DGs 2\&3\\
\midrule
Rated Active Power & $P_i^*$ & 1400 W & 700 W\\
Rated Reactive Power & $Q_i^*$ & 800 VAr & 400 VAr \\
$P-\omega$ Droop Coeff. & $m_i$ & $2.5 \cdot 10^{-3}$ $\frac{\rm rad}{\rm{W\cdot s}}$ & $5 \cdot 10^{-3}$ $\frac{\rm rad}{\rm{W\cdot s}}$\\
$Q-E$ Droop Coeff. & $n_i$ & $1.5 \cdot 10^{-3}$ $\frac{\rm V}{\rm{VAr}}$ & $3 \cdot 10^{-3}$ $\frac{\rm V}{\rm{VAr}}$\\
Int. Frequency Gain  & $k_i$ & 1.7\,s & 1.7\,s\\
Int. Voltage Gain & $\kappa_i$ & 1\,s & 1\,s\\
\bottomrule
\end{tabular}
}
\label{Tab:Parameters}
\end{center}
\end{table}

\subsection{Study 1: Controller Performance}
\label{Sec:Study1}


\begin{figure*}[!ht]
        \centering
        \begin{subfigure}[!ht]{0.3\textwidth}
                \includegraphics[width=\columnwidth]{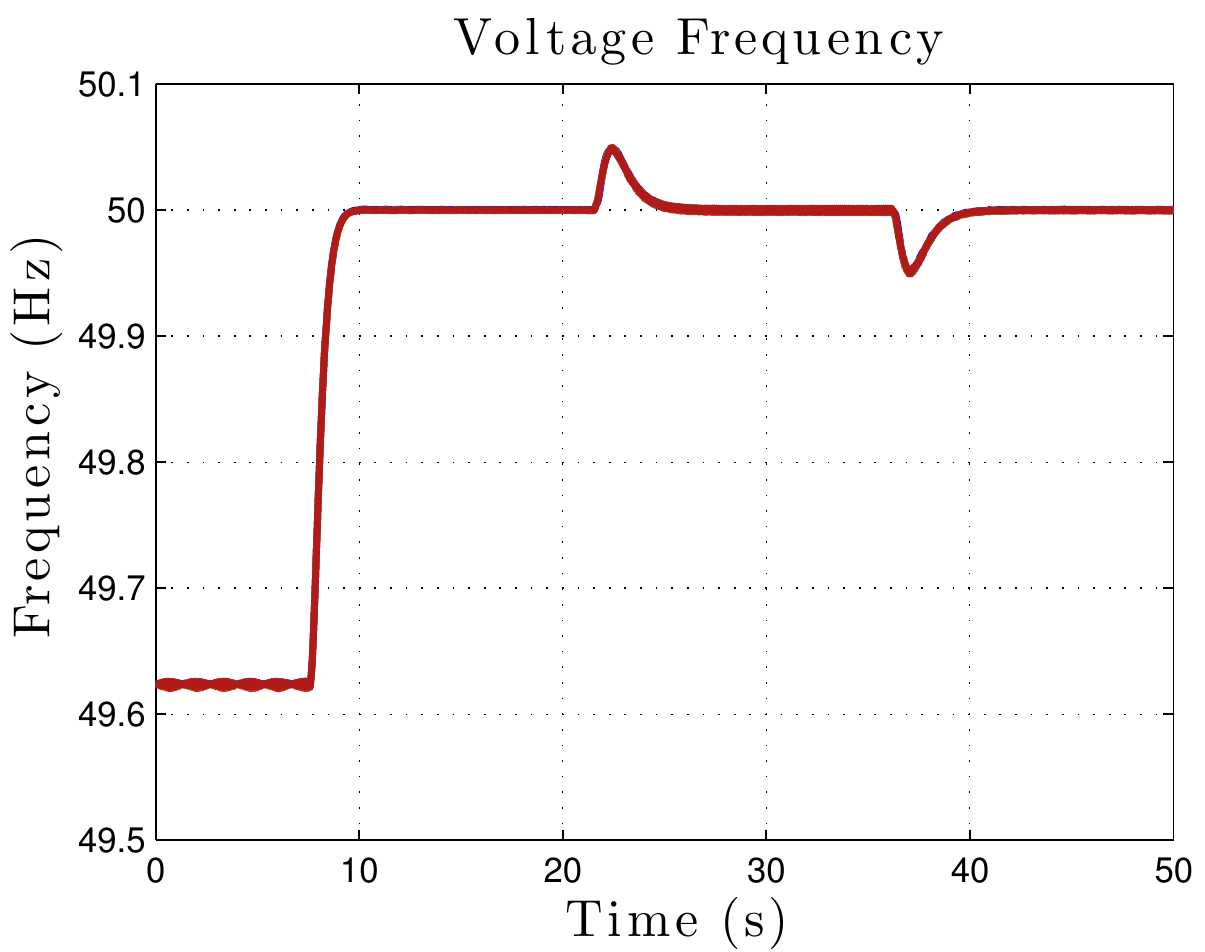}
                \label{Fig:1af}
        \end{subfigure}~
        \begin{subfigure}[!ht]{0.3\textwidth}
                \includegraphics[width=\columnwidth]{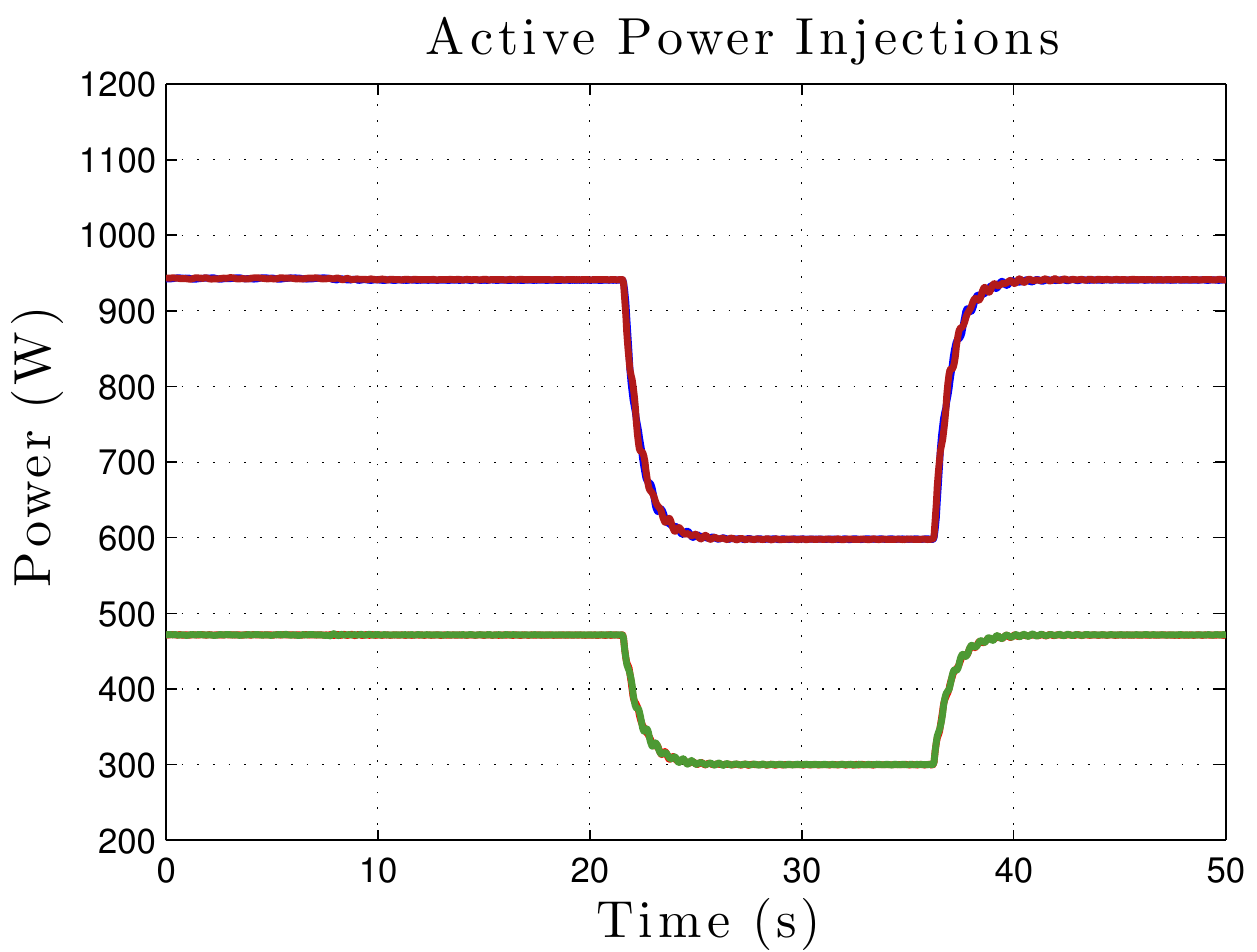}
                \label{Fig:1aP}
        \end{subfigure}~
        \begin{subfigure}[!ht]{0.3\textwidth}
                \includegraphics[width=\columnwidth]{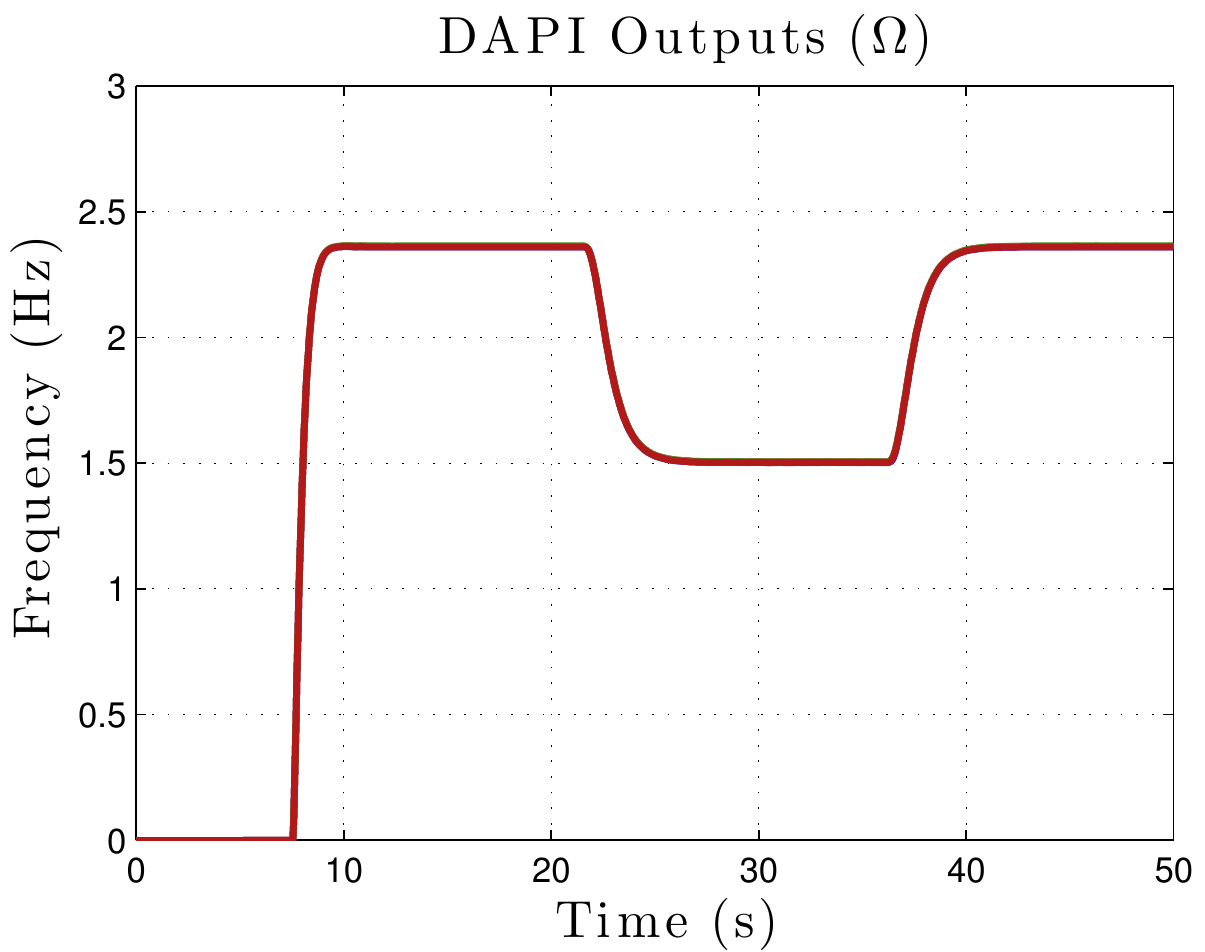}
                \label{Fig:1aO}
        \end{subfigure}
        
        \begin{subfigure}[!ht]{0.3\textwidth}
                \includegraphics[width=\columnwidth]{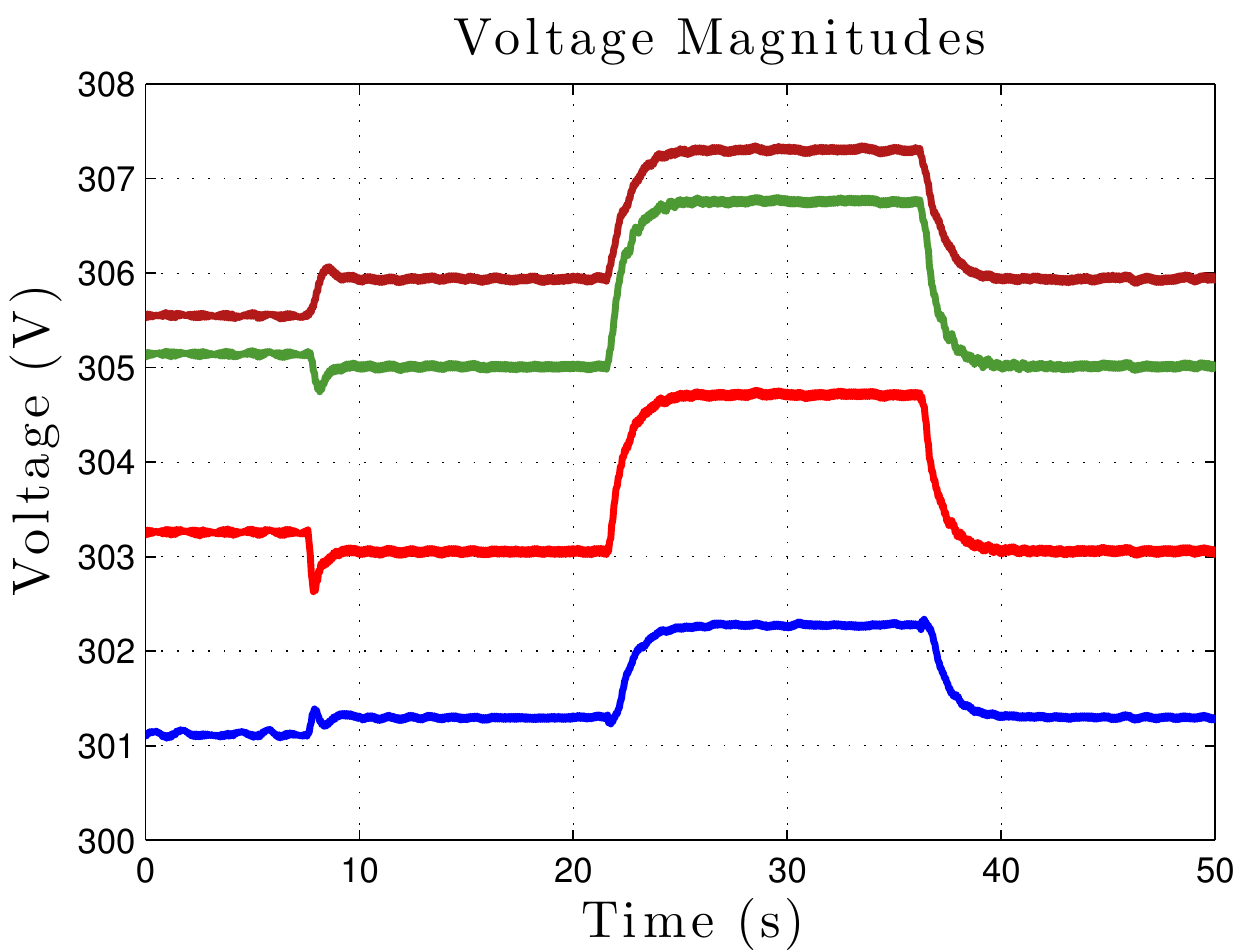}
                \label{Fig:1aE}
        \end{subfigure}~
         \begin{subfigure}[!ht]{0.3\textwidth}
                \includegraphics[width=\columnwidth]{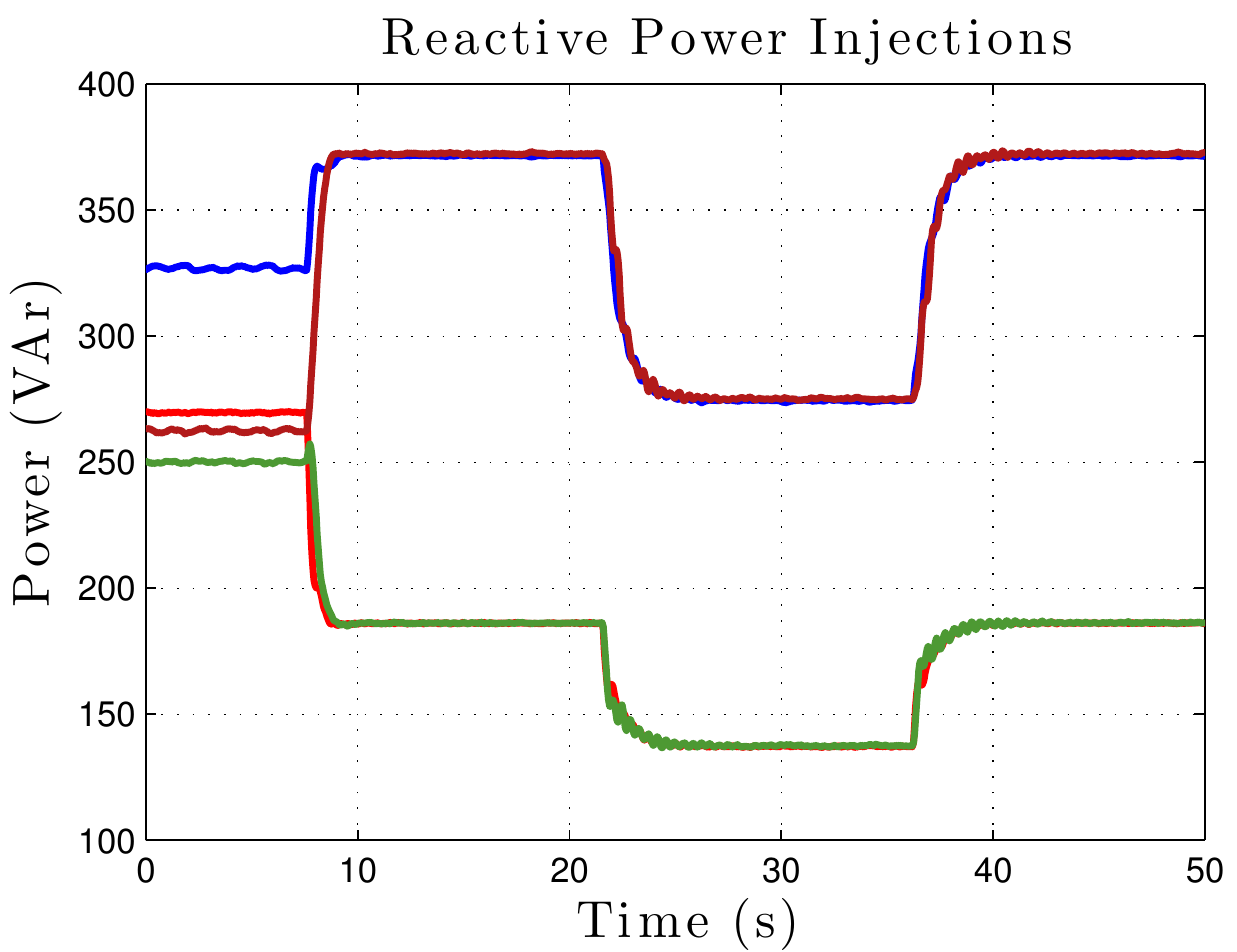}
                \label{Fig:1aQ}
        \end{subfigure}~
        \begin{subfigure}[!ht]{0.3\textwidth}
                \includegraphics[width=\columnwidth]{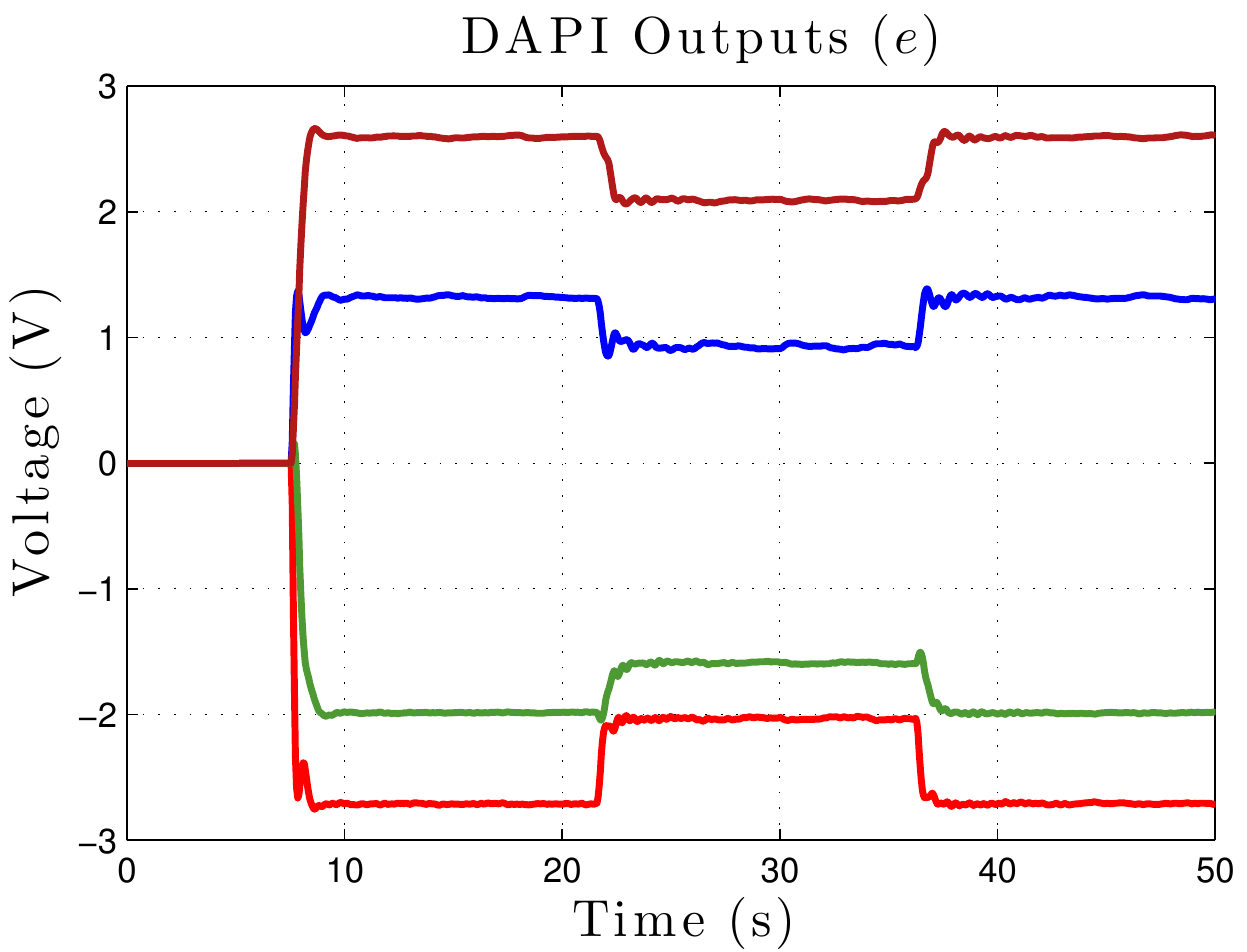}
                \label{Fig:1ase}
        \end{subfigure}
        \captionsetup{justification=raggedright,singlelinecheck=false}
        \caption{Study 1a -- reactive power sharing without voltage regulation, with control parameters $b = 50$\,V, $\beta_1 = \beta_2 = \beta_3 = \beta_4 = 0$. In correspondence with Figure \ref{Fig:VoltCharSharing} of Section \ref{Sec:ReactiveSharing}, the quality of voltage regulation is quite poor.}\label{Fig:1a}
\end{figure*}


\begin{figure*}
        \centering
        \begin{subfigure}[!ht]{0.3\textwidth}
                \includegraphics[width=\columnwidth]{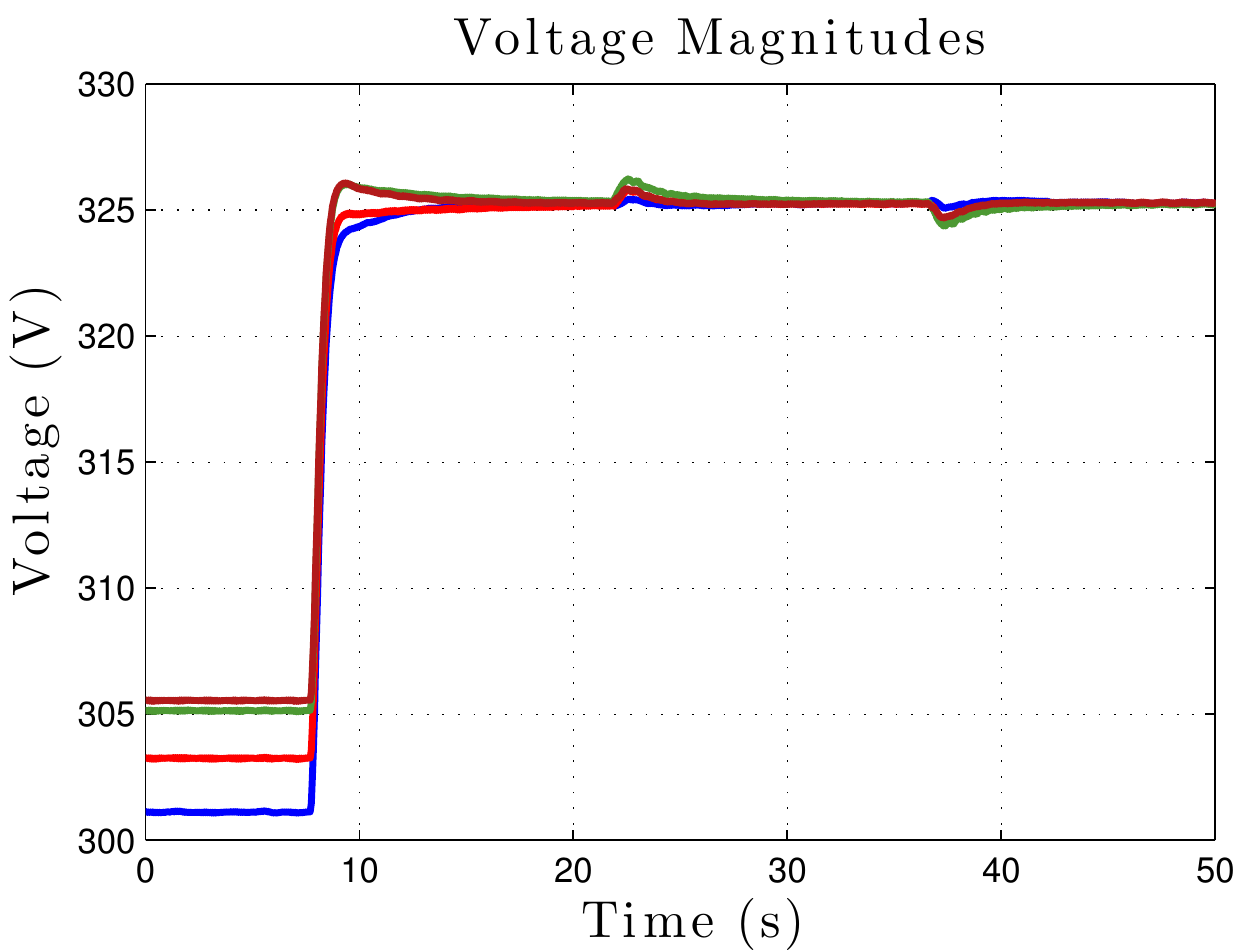}
                \label{Fig:1bE}
        \end{subfigure}~
        \begin{subfigure}[!ht]{0.3\textwidth}
                \includegraphics[width=\columnwidth]{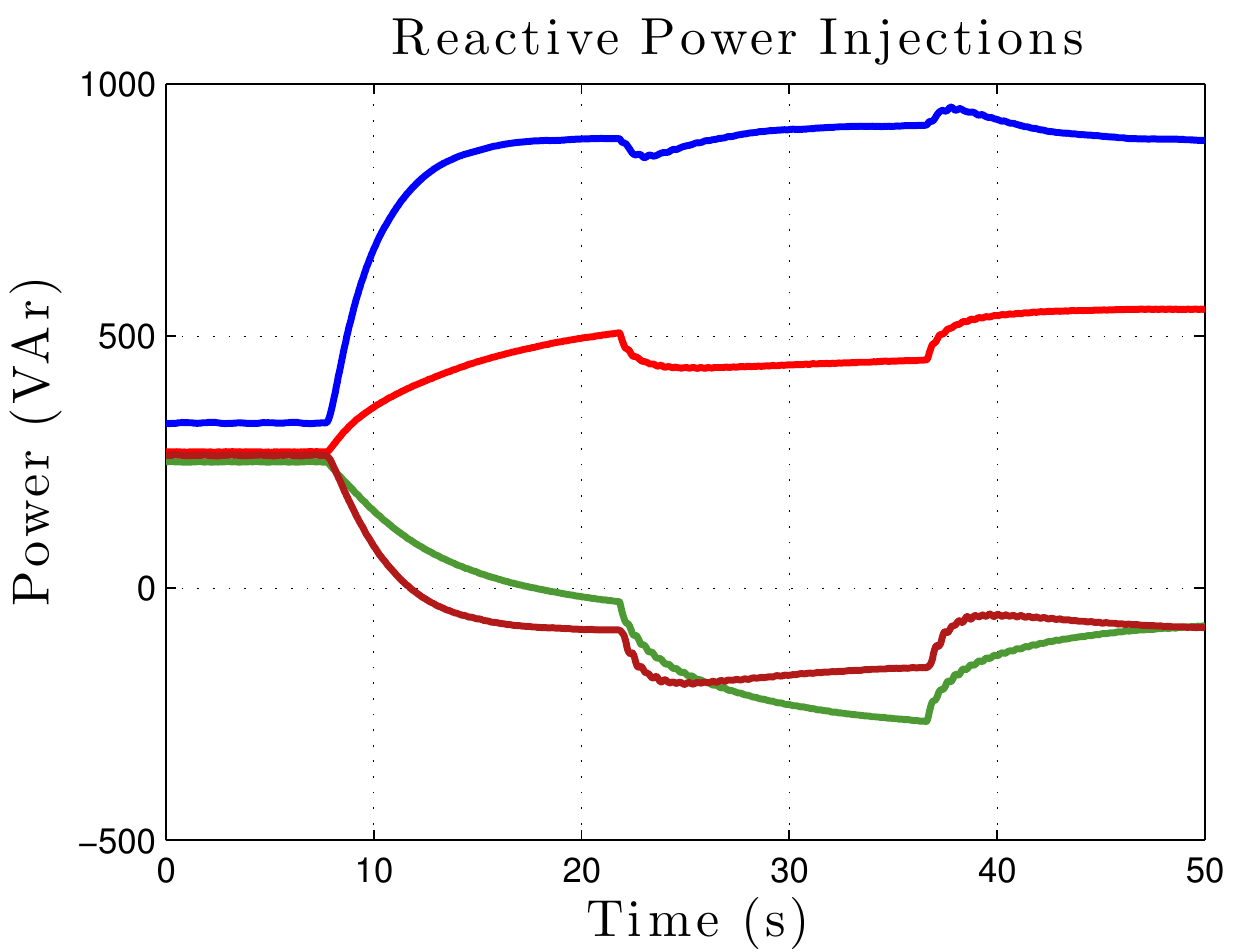}
                \label{Fig:1bQ}
        \end{subfigure}~
        \begin{subfigure}[!ht]{0.3\textwidth}
                \includegraphics[width=\columnwidth]{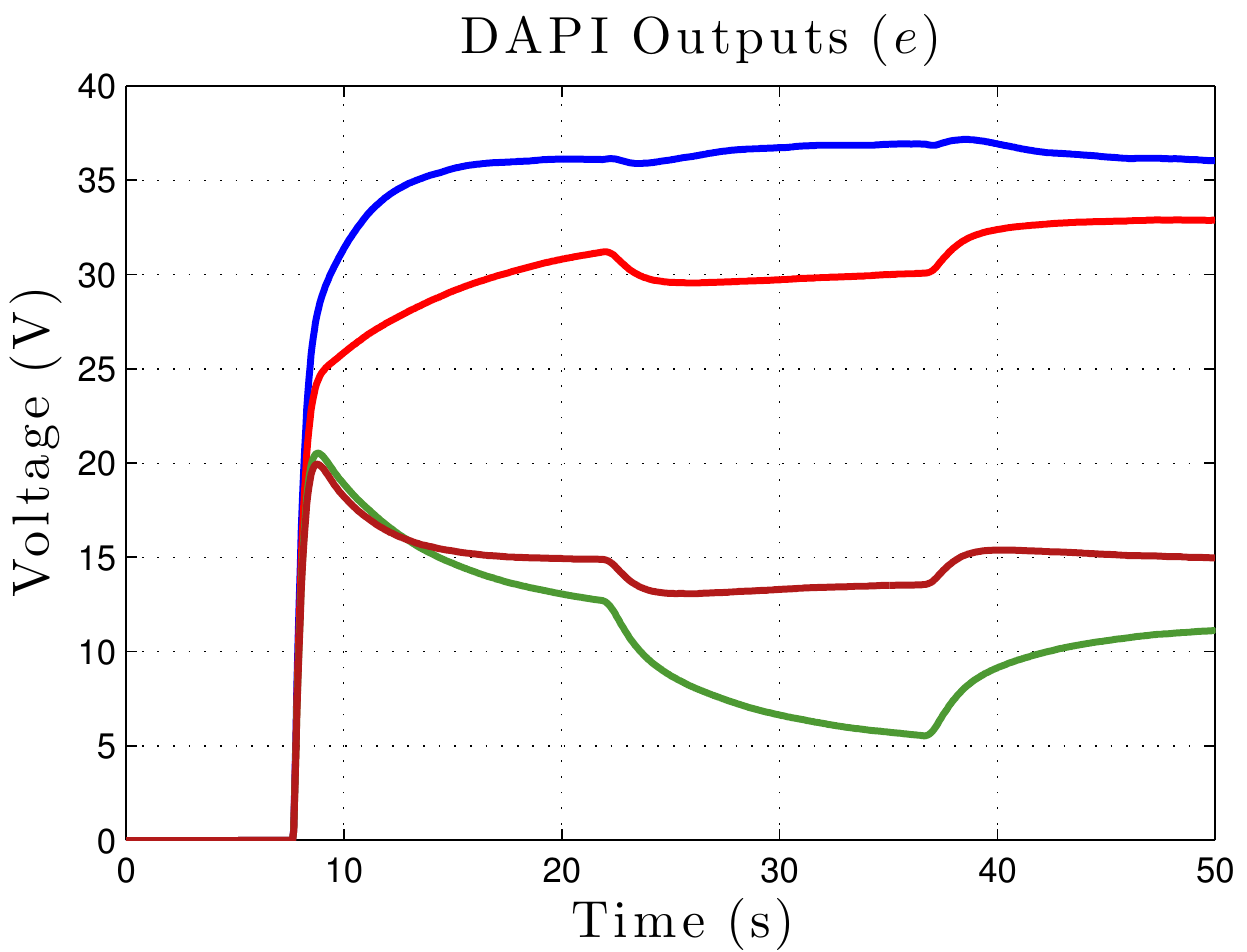}
                \label{Fig:1bse}
        \end{subfigure}
        \captionsetup{justification=raggedright,singlelinecheck=false}
        \caption{Study 1b -- voltage regulation without reactive power sharing, with parameters $b = 0$\,V, $\beta_1 = \beta_2 = \beta_3 = \beta_4 = 2.2$. In correspondence with Figure \ref{Fig:VoltChar} of Section \ref{Sec:ReactiveSharing}, the quality of reactive power sharing is quite poor.}\label{Fig:1b}
\end{figure*}


\begin{figure*}[t!]
        \centering
%
        \begin{subfigure}[!ht]{0.3\textwidth}
                \includegraphics[width=\columnwidth]{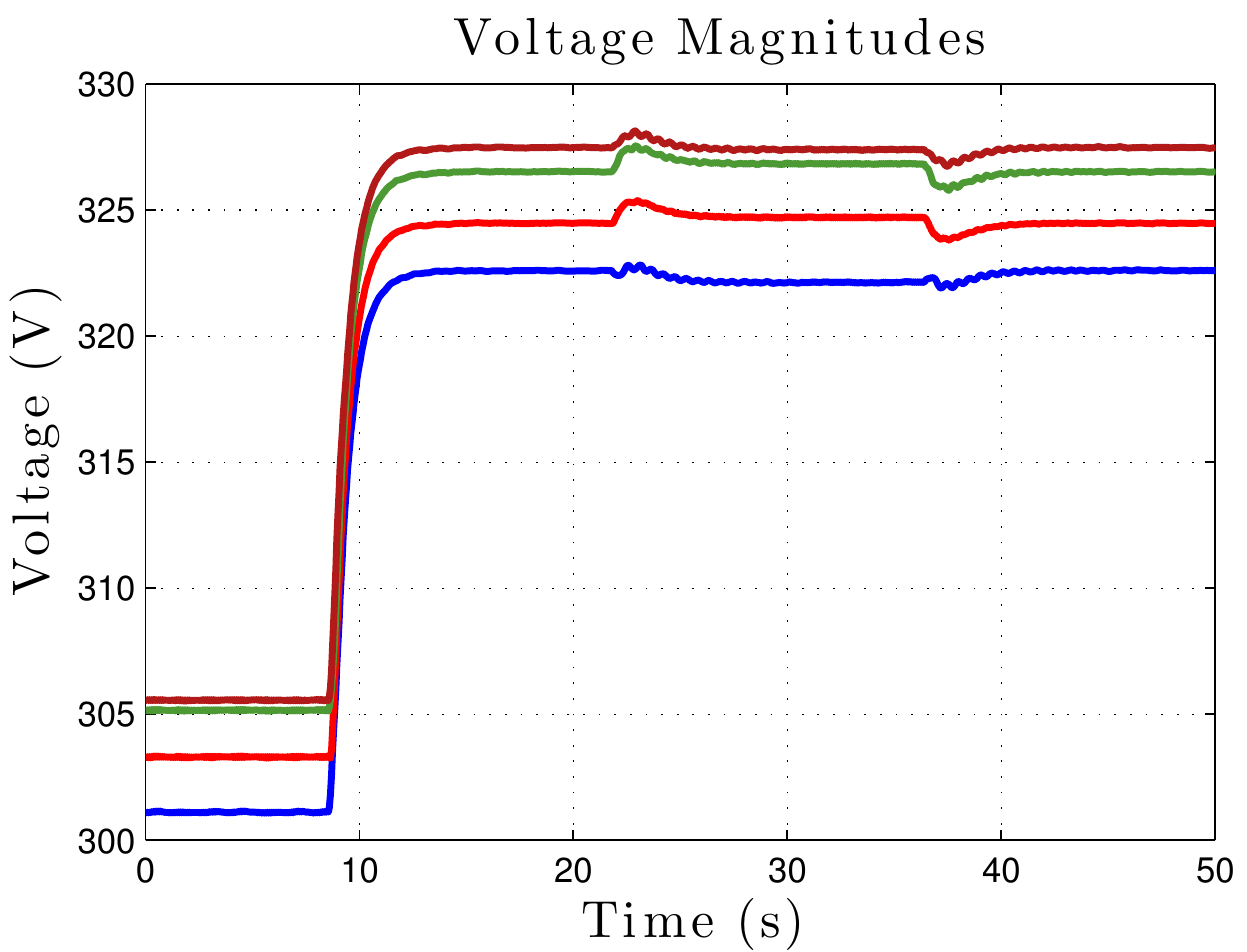}
                \label{Fig:1cE}
        \end{subfigure}
        ~ 
        \begin{subfigure}[!ht]{0.3\textwidth}
                \includegraphics[width=\columnwidth]{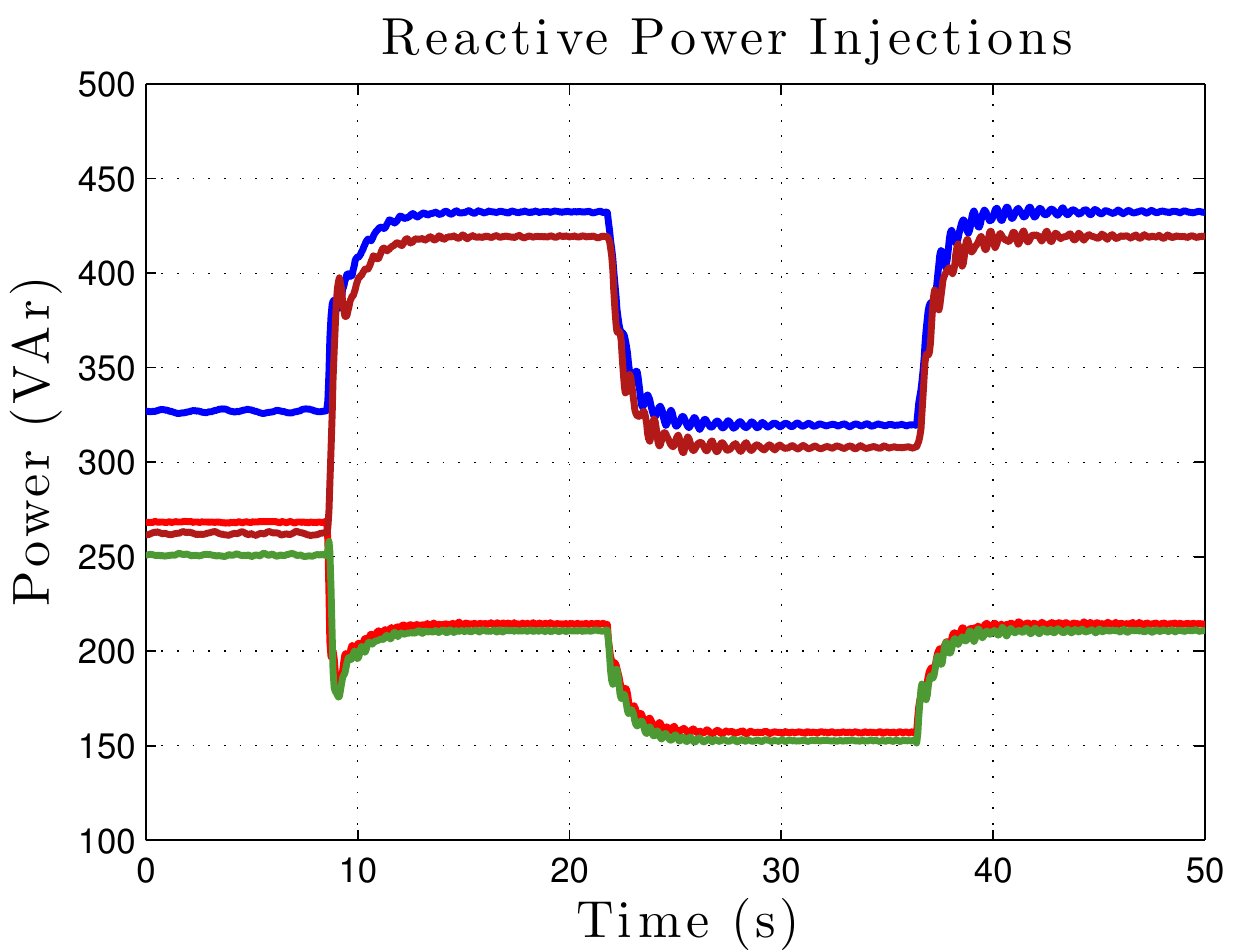}
                \label{Fig:1cQ}
        \end{subfigure}
        ~
        \begin{subfigure}[!ht]{0.3\textwidth}
                \includegraphics[width=\columnwidth]{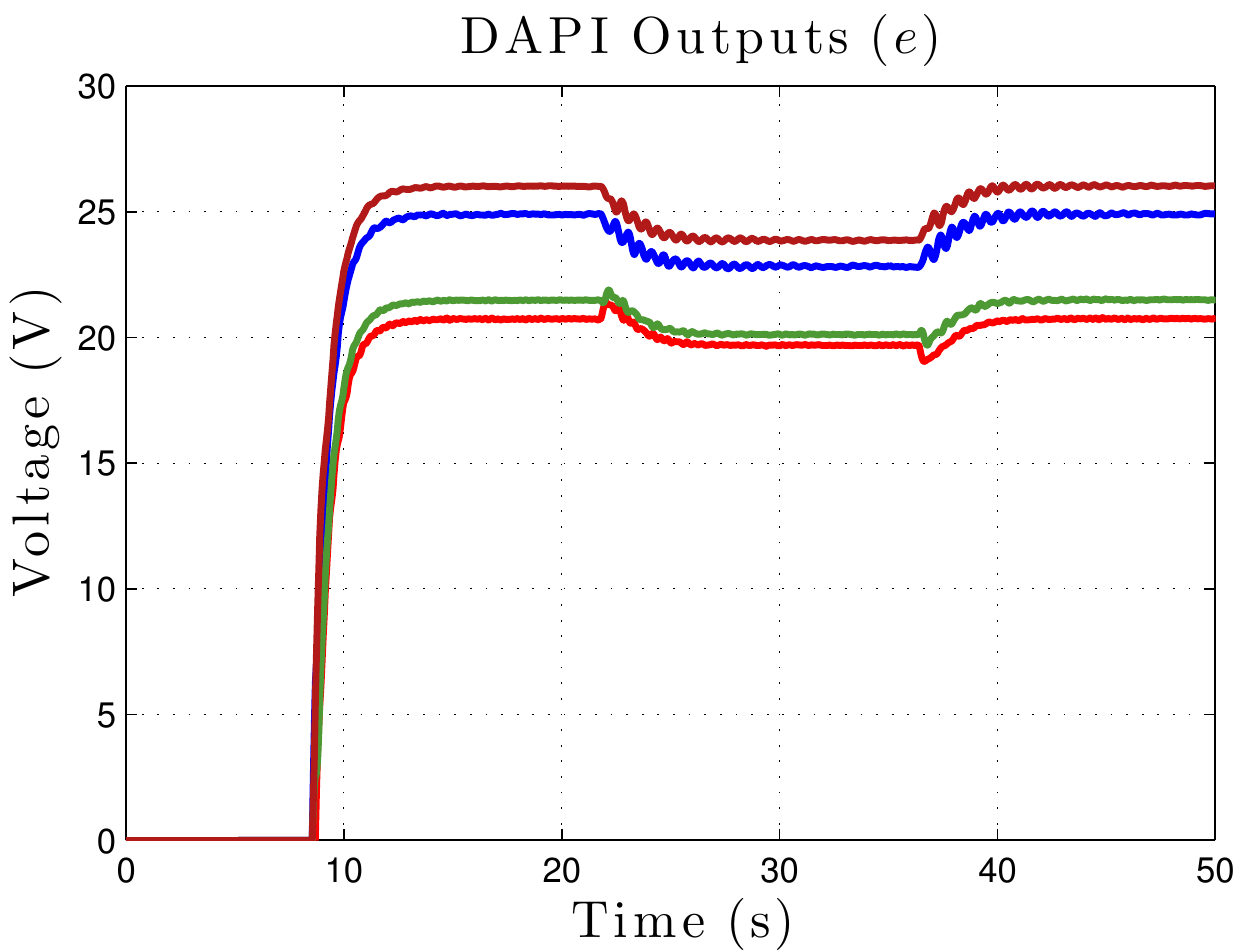}
                \label{Fig:1cse}
        \end{subfigure}
        \captionsetup{justification=raggedright,singlelinecheck=false}
        \caption{Study 1c -- a compromise between voltage regulation and reactive power sharing, with control parameters $\beta_1 = \beta_2 = \beta_3 = \beta_4 =  1.2$, $b = 180$\,V.}\label{Fig:1c}
\end{figure*}


\begin{figure*}
        \centering
%
        \begin{subfigure}[!ht]{0.3\textwidth}
                \includegraphics[width=\columnwidth]{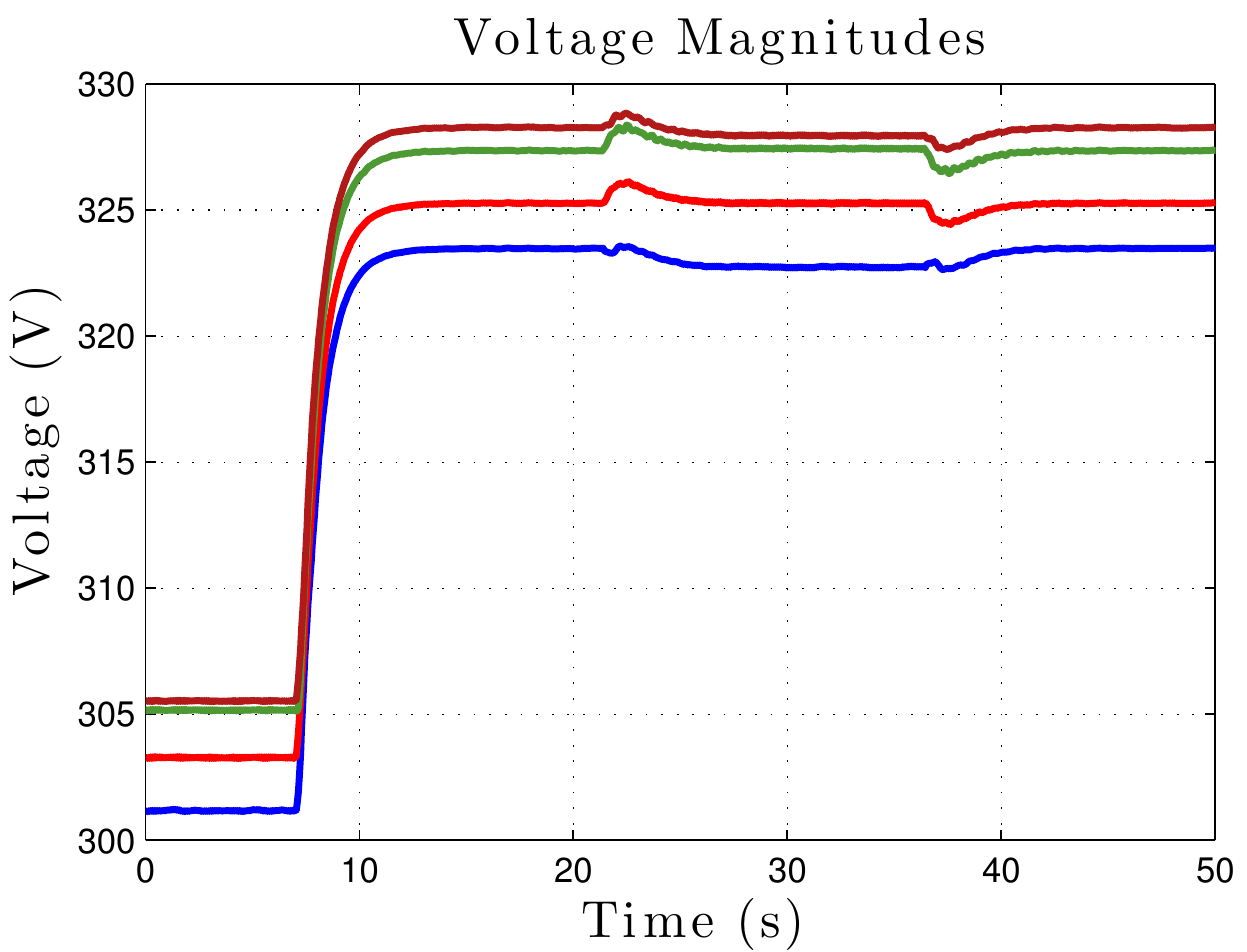}
                \label{Fig:1eE}
        \end{subfigure}
        ~ 
        \begin{subfigure}[!ht]{0.3\textwidth}
                \includegraphics[width=\columnwidth]{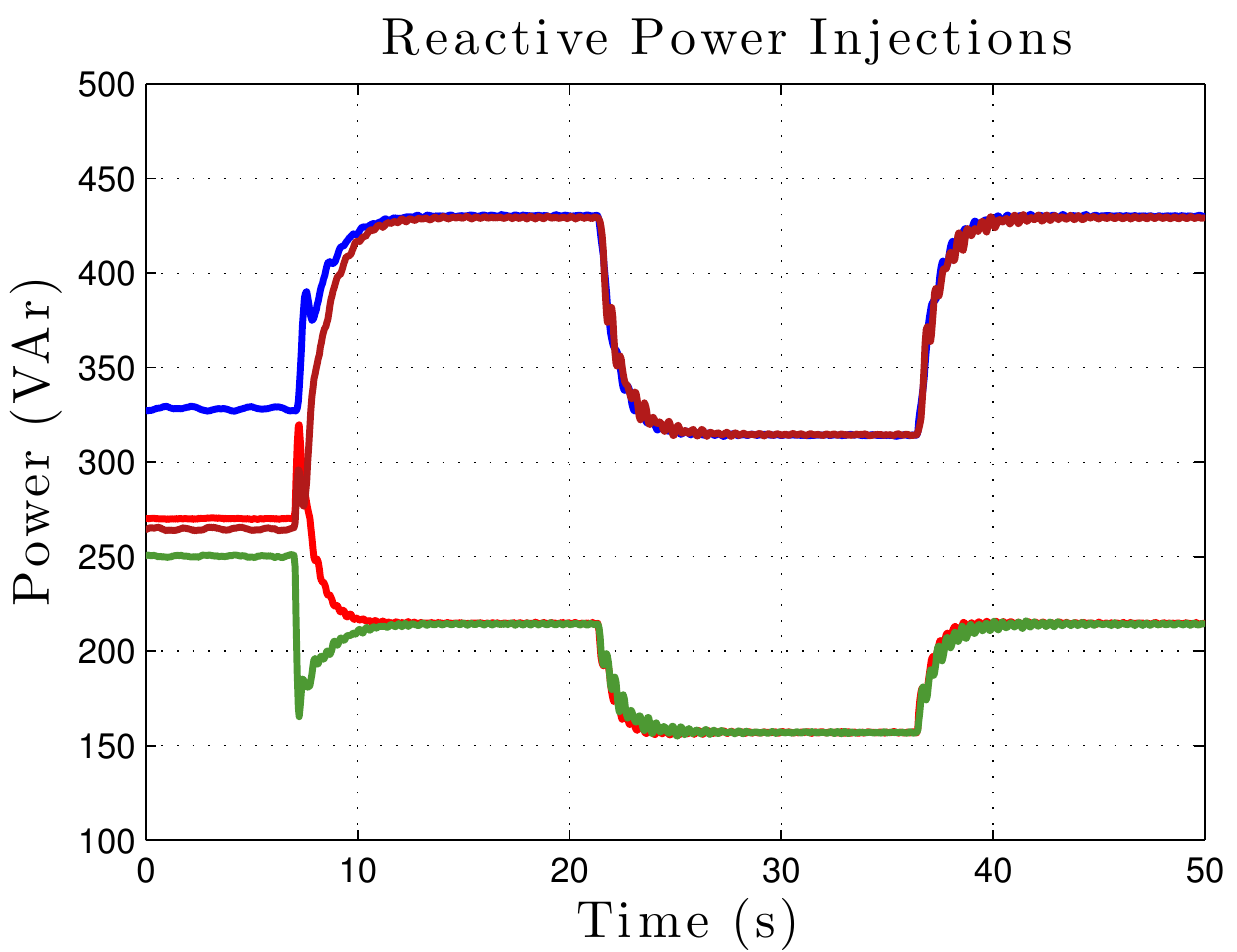}
                \label{Fig:1eQ}
        \end{subfigure}
        ~
        \begin{subfigure}[!ht]{0.3\textwidth}
                \includegraphics[width=\columnwidth]{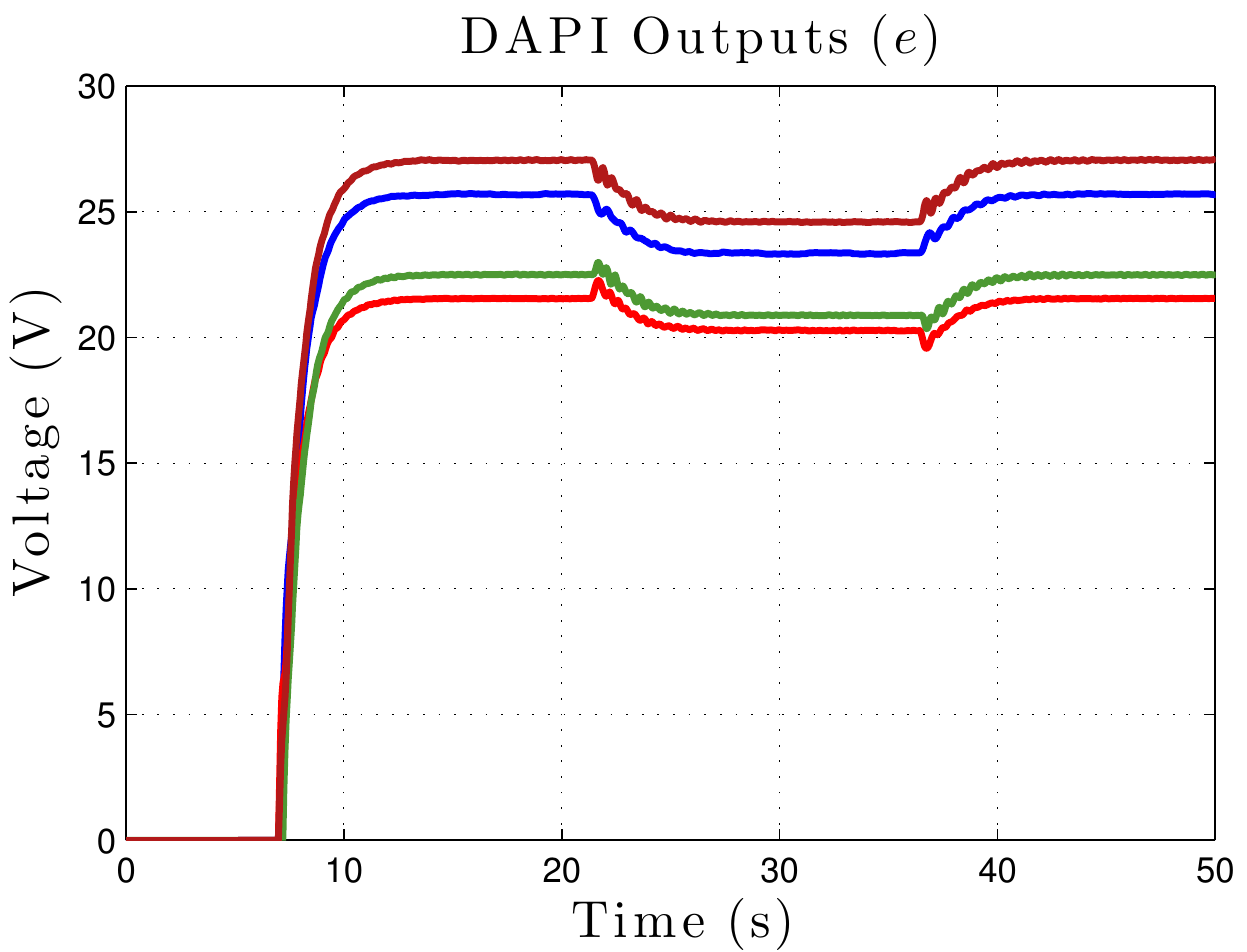}
                \label{Fig:1ese}
        \end{subfigure}
        \captionsetup{justification=raggedright,singlelinecheck=false}
        \caption{Study 1d -- accurate reactive power sharing and good voltage regulation, with control parameters $\beta_1 = \beta_3 = \beta_4 = 0$, $\beta_2 = 4$ and $b = 100$\,V.}\label{Fig:1d}
\end{figure*}

Studies in this section illustrate the performance of the DAPI controllers \eqref{Eq:Primary}--\eqref{Eq:SecondaryReactive} under various controller tunings. In all four sub-studies 1a--1d, only primary droop controllers \eqref{Eq:Primary} and \eqref{Eq:PrimaryReactive} are running up to $t = 7$\,s, at which time the secondary controllers \eqref{Eq:Secondary} and \eqref{Eq:SecondaryReactive} are activated. The local load at DG unit 4 is detached at $t=22$\,s, then reattached at $t = 36$\,s.

First considering the frequency dynamics in the top portion of Figure \ref{Fig:1a}, the frequency deviation experienced under primary droop control is quickly eliminated by the DAPI controller \eqref{Eq:Primary}--\eqref{Eq:Secondary}, and frequency regulation is maintained throughout load changes with minimal transients. Active power is accurately shared amongs the heterogeneous DGs throughout the entire runtime. This robust frequency and active power behavior is identical in all other sub-studies, and thus we omit the plots due to space considerations. The remainder of studies 1a and 1b in Figures \ref{Fig:1a} and \ref{Fig:1b} highlight the conclusions drawn in Section \ref{Sec:ReactiveSharing} regarding the limitations of voltage secondary control. 

\paragraph{Tuning for pure reactive power sharing}
Figure \ref{Fig:1a} shows results for the voltage DAPI controller \eqref{Eq:PrimaryReactive}--\eqref{Eq:SecondaryReactive} tuned for power sharing ($\beta_i = 0$, $b = 50$\,V), with no attempt at voltage regulation, as in Case 1 of Section \ref{Section:QDAPI}. While reactive power is shared accurately, voltage magnitudes deviate from their nominal values $E^* = 325.3$\,V (cf. Figure \ref{Fig:VoltCharSharing}).
\paragraph{Tuning for pure voltage regulation}
Conversely, Figure \ref{Fig:1b} reports results for the same controller tuned to regulate voltage levels without reactive power sharing ($\beta_i = 2.2$, $b = 0$\,V), as in Case 2 of Section \ref{Section:QDAPI}. While voltage levels are tightly regulated, reactive power sharing among the units is poor (cf. Figure \ref{Fig:VoltChar}). As explained in Section \ref{Sec:ReactiveSharing}, the poor performance in Figure \ref{Fig:1b} is a general property of \emph{all} voltage controllers strategies that exactly regulate DG voltage levels.
\paragraph{Compromised tuning}
Figure \ref{Fig:1c} displays the results for Study 1c, in which the DAPI controllers \eqref{Eq:Primary}--\eqref{Eq:SecondaryReactive} are implemented with $b = 180$\,V and uniform controller gains $\beta_i = 1.2$, as in Case 3 of Section \ref{Section:QDAPI}. Considering the voltage dynamics, the voltage DAPI controller \eqref{Eq:PrimaryReactive}--\eqref{Eq:SecondaryReactive} achieves a \emph{compromise} between voltage regulation and reactive power sharing. Voltage magnitudes are roughly clustered around $E^* = 325.3$\,V, while reactive power is approximately shared.
\paragraph{Smart tuning}
Figure \ref{Fig:1d} displays the results for Study 1d, in which the DAPI controllers \eqref{Eq:Primary}--\eqref{Eq:SecondaryReactive} are implemented with $b = 100$\,V and $\beta_2 = 4$, $\beta_1 = \beta_3 = \beta_4 = 0$, in accordance with the discussion of Case 4 in Section \ref{Section:QDAPI}. In comparison with the voltage dynamics of Study 1c, the voltage regulation in Figure \ref{Fig:1d} shows a slight improvement, while the reactive power sharing is noticeably improved, maintaining accurate sharing through load changes and during transients. Note that this performance improvement has been achieved while \emph{reducing} the controller gain $b$ which enforces reactive power sharing. Due to this reduction in controller gain, the ringing in the reactive power signal during transients is noticeably improved from Study 1c to Study 1d, { in agreement with the stability and root locus analyses of Section \ref{Sec:Stability}}.

\subsection{Study 2: Communication Link Failure}
\label{Sec:Study2}

\begin{figure}[!t]
        \centering
        \begin{subfigure}[!ht]{0.48\columnwidth}
                \includegraphics[width=0.9\columnwidth]{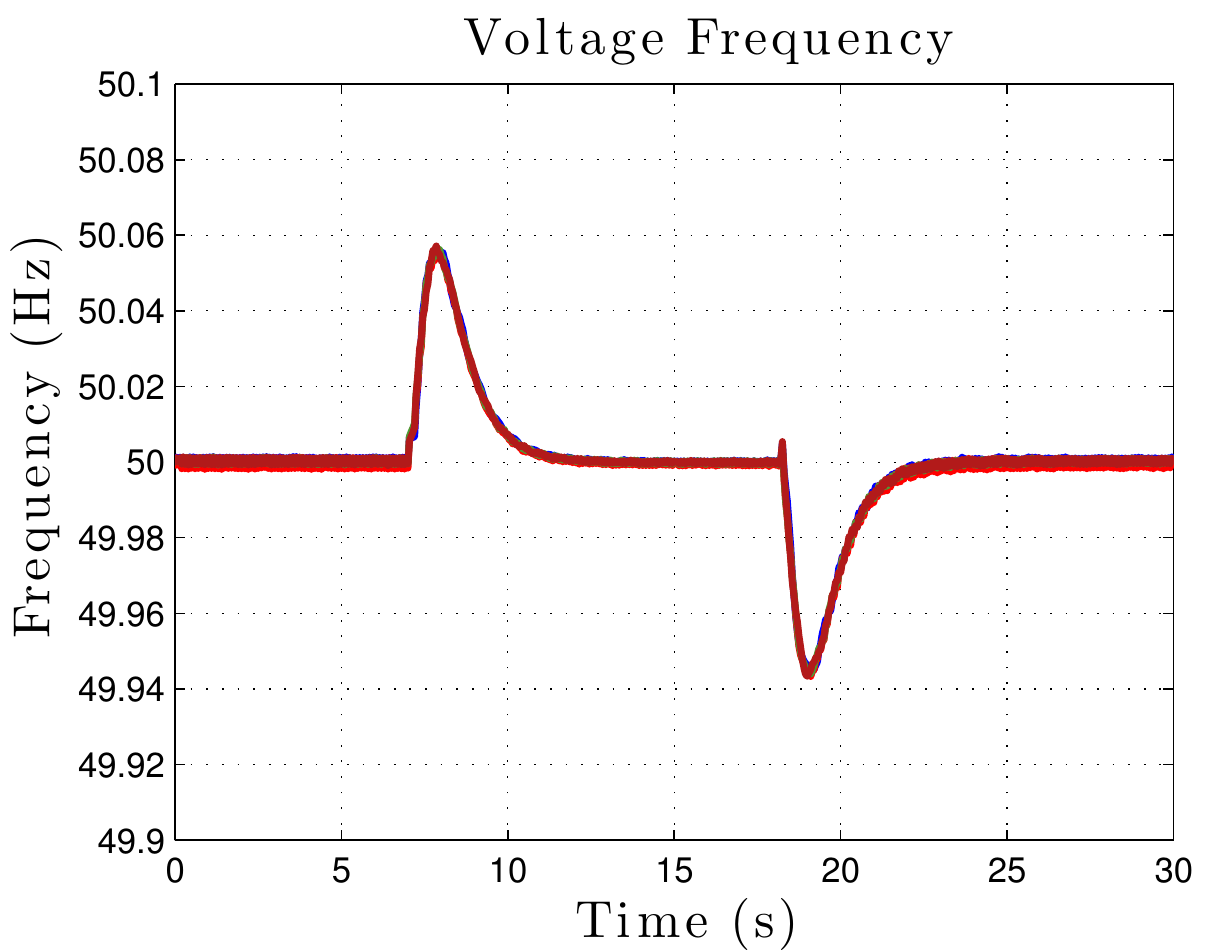}
                \label{Fig:2f}
        \end{subfigure}~
        \begin{subfigure}[!ht]{0.48\columnwidth}
                \includegraphics[width=0.9\columnwidth]{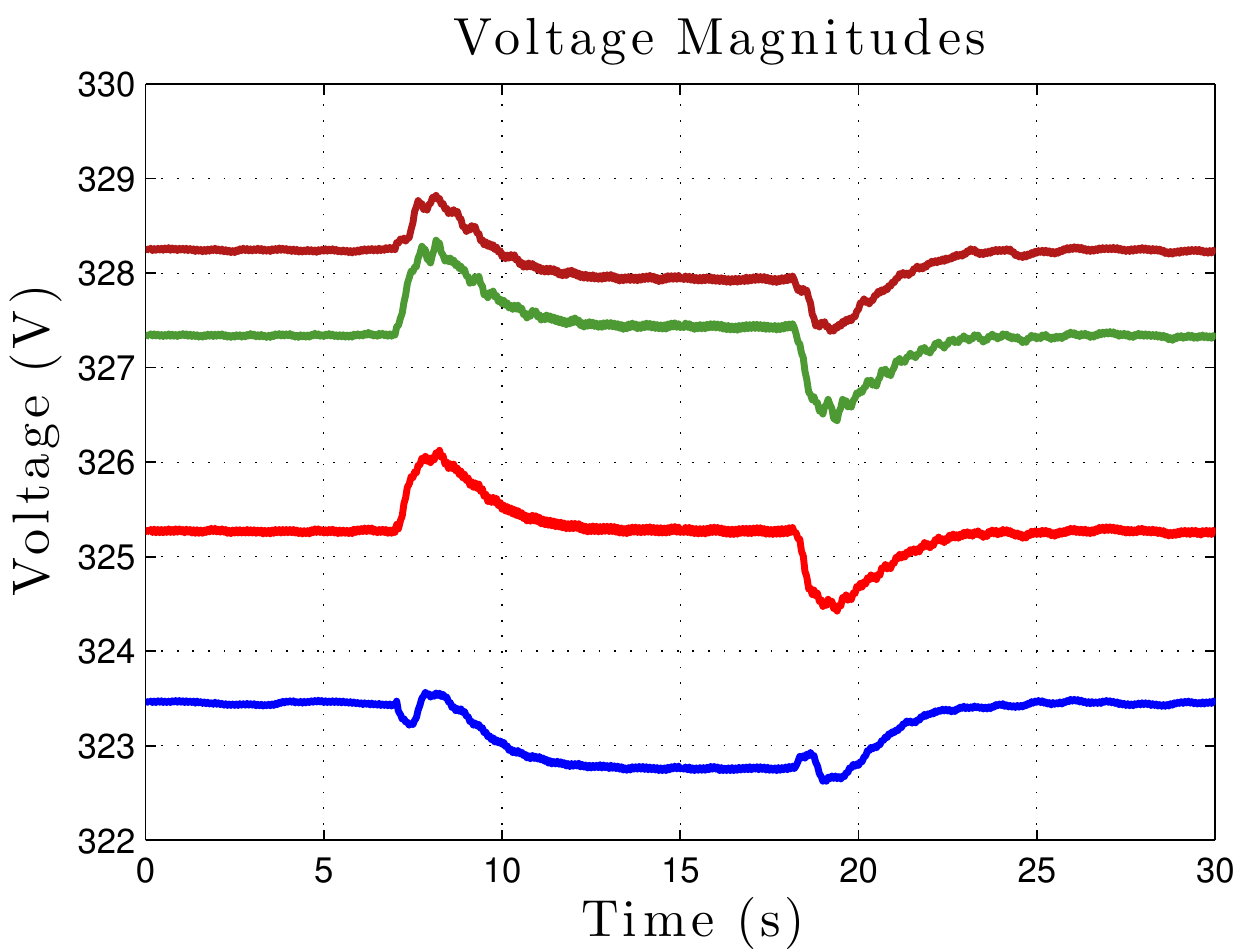}
                \label{Fig:2E}
        \end{subfigure}
        
        \begin{subfigure}[!ht]{0.48\columnwidth}
                \includegraphics[width=0.9\columnwidth]{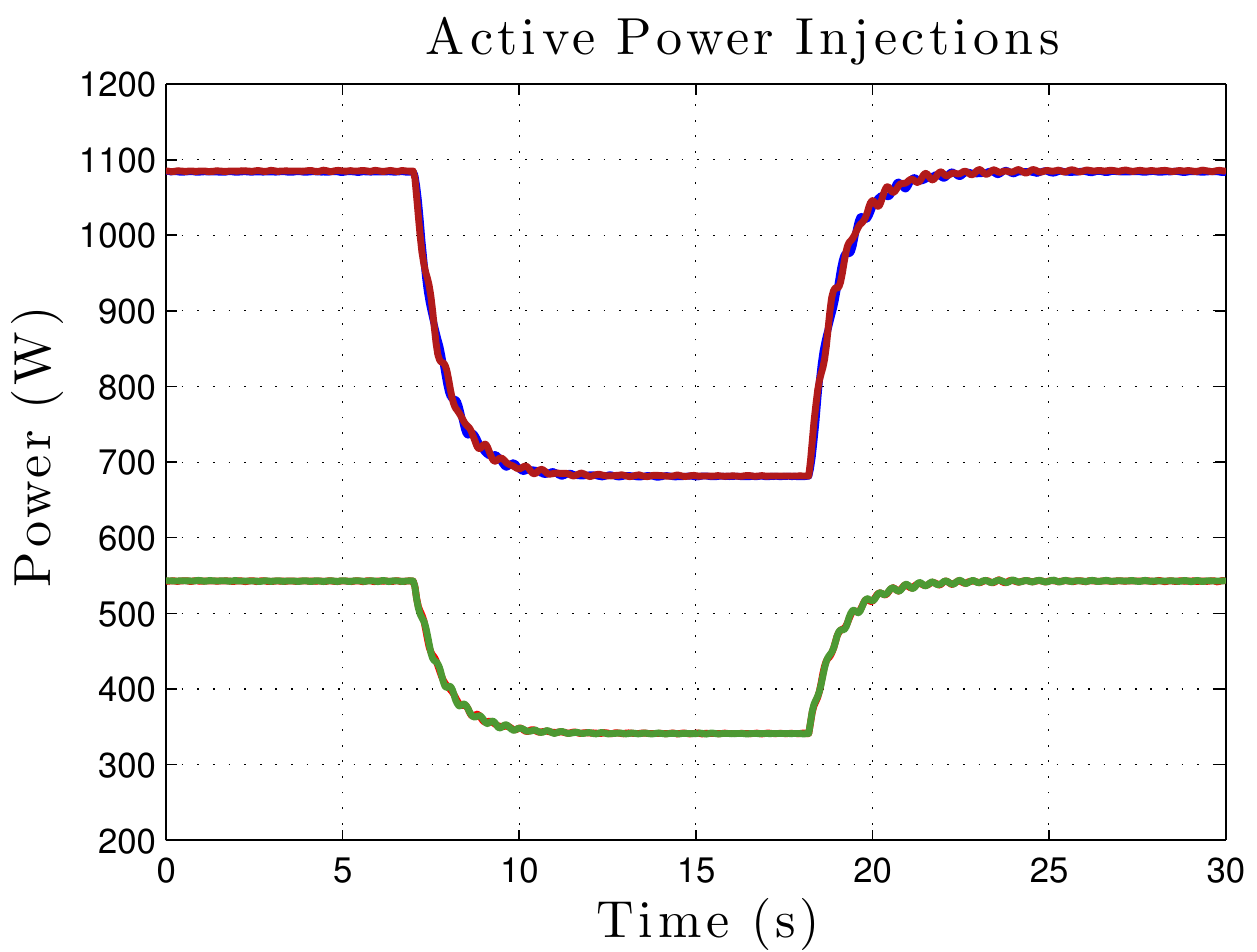}
                \label{Fig:2P}
        \end{subfigure}~
        \begin{subfigure}[!ht]{0.48\columnwidth}
                \includegraphics[width=0.9\columnwidth]{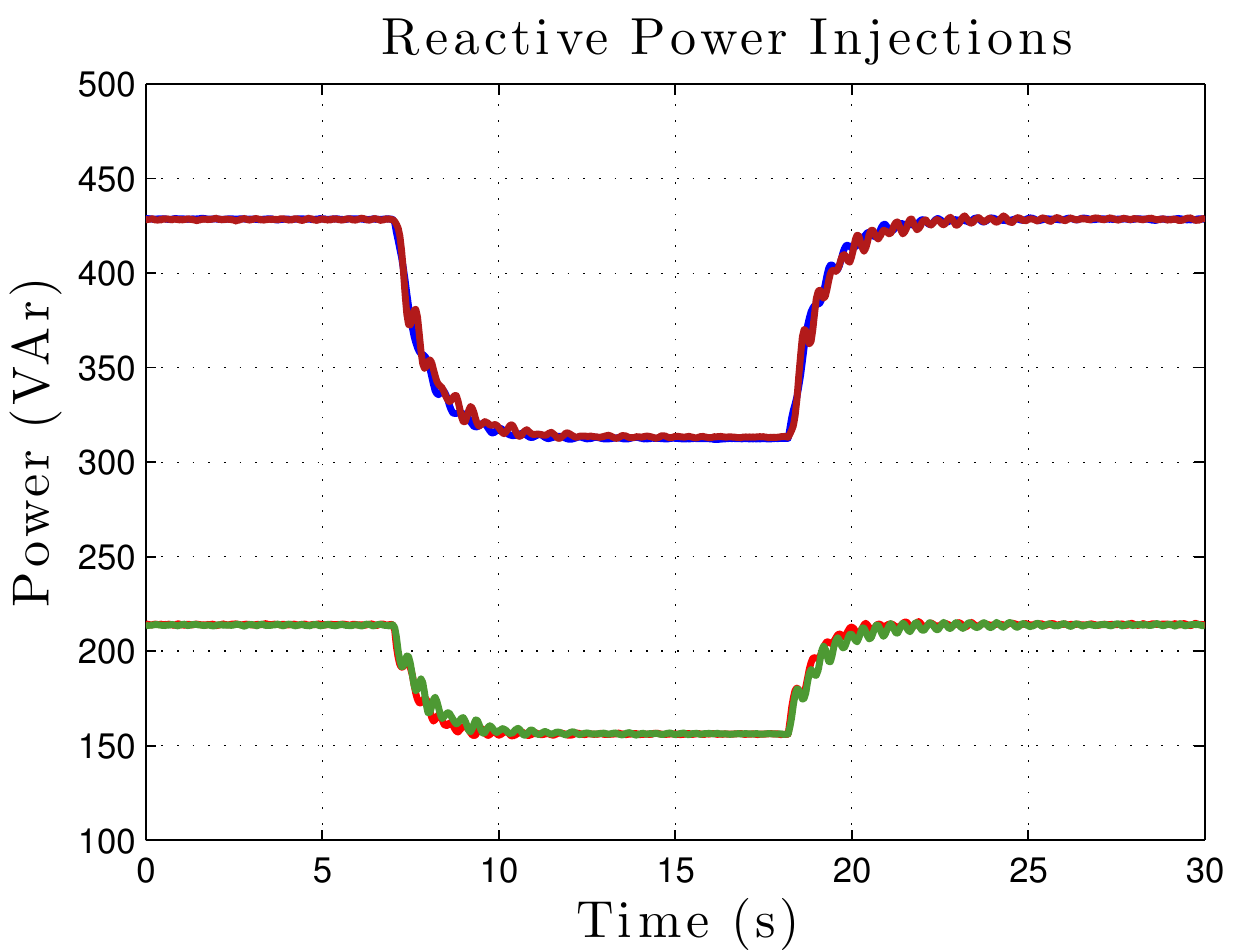}
                \label{Fig:2Q}
        \end{subfigure}
        
        \begin{subfigure}[!ht]{0.48\columnwidth}
                \includegraphics[width=0.9\columnwidth]{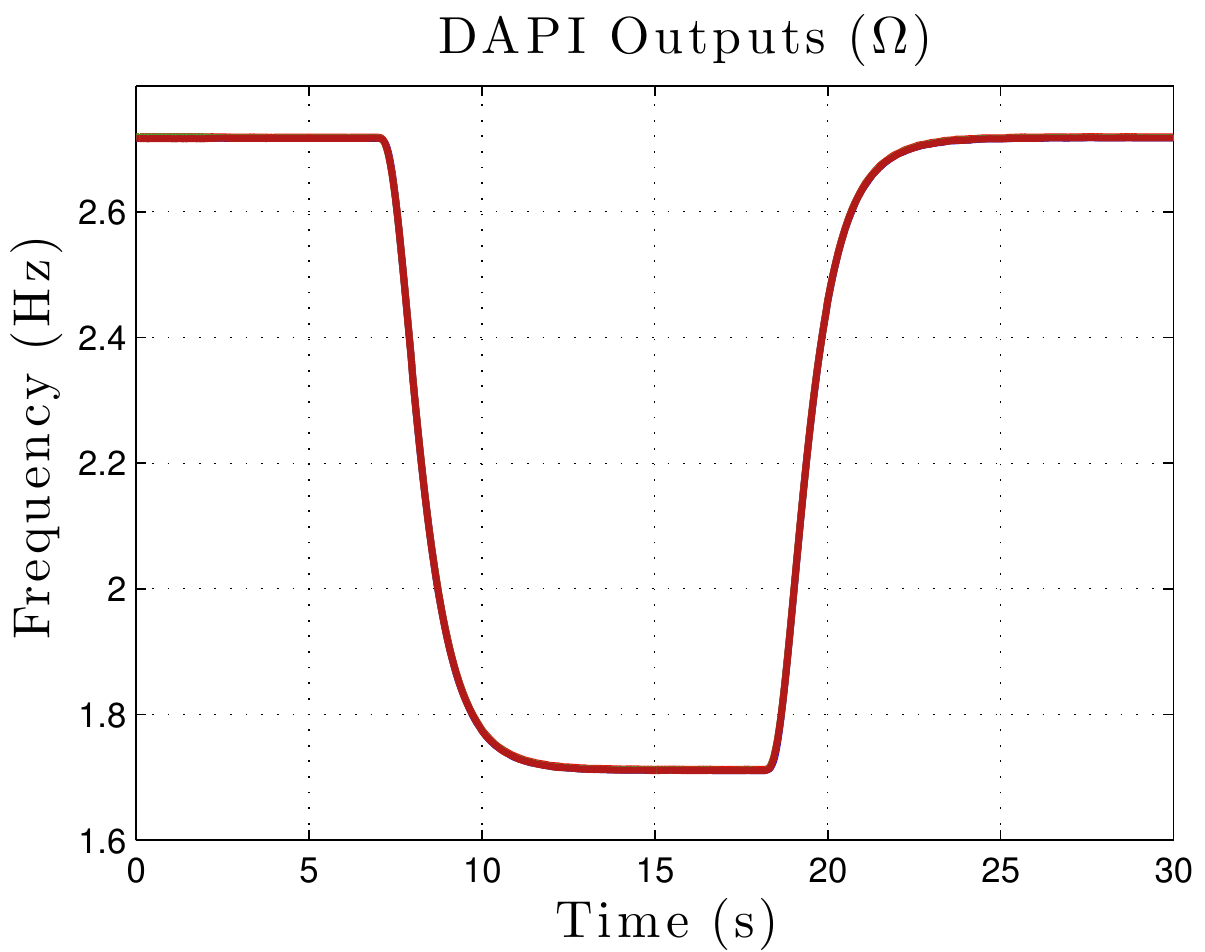}
                \label{Fig:2O}
        \end{subfigure}~
        \begin{subfigure}[!ht]{0.48\columnwidth}
                \includegraphics[width=0.9\columnwidth]{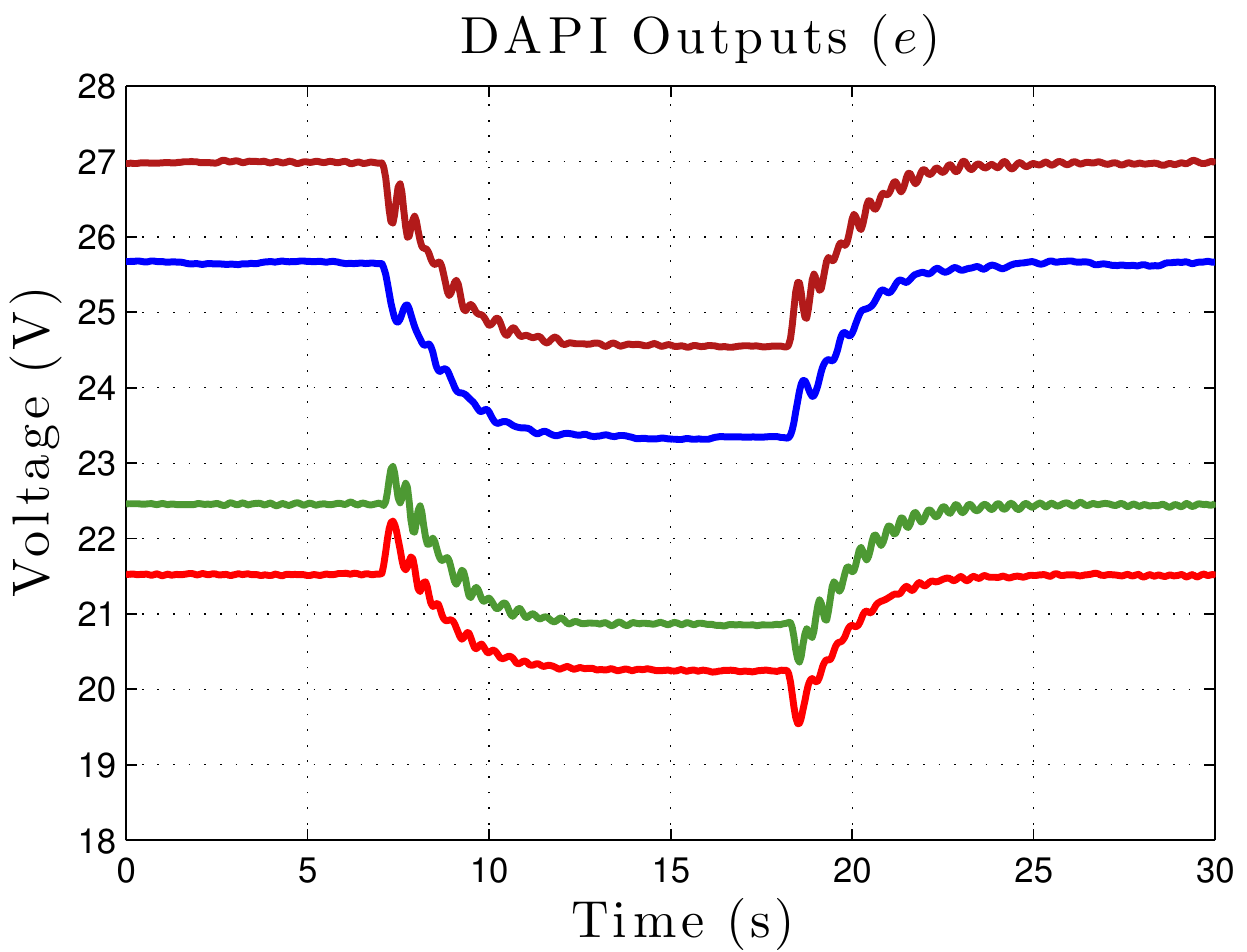}
                \label{Fig:2se}
        \end{subfigure}
        \captionsetup{justification=raggedright,singlelinecheck=false}
        \caption{Study 2 -- DAPI performance under communication link failure. Control parameters are the same as in Study 1d.}\label{Fig:2}
\end{figure}

In this study the communication link (Figure \ref{Fig:Schematic}) between DG units 3 and 4 fails at $t = 2$\,s. At $t = 7$\,s the local load at unit 4 is detached, and is reattached at $t = 18$\,s. Control parameters are the same as in Study 1d. As the results in Figure \ref{Fig:2} show, the DAPI controllers \eqref{Eq:Primary}--\eqref{Eq:SecondaryReactive} maintain the high performance from Study 1d despite the absence of the communication link between DG units 3 and 4 (cf. Remark \ref{Rem:Interp}).

\subsection{Study 3: Non-Uniform Controller Gains}
\label{Sec:Study5}

\begin{figure}[!ht]
        \centering
        \begin{subfigure}[!ht]{0.47\columnwidth}
                \includegraphics[width=\columnwidth]{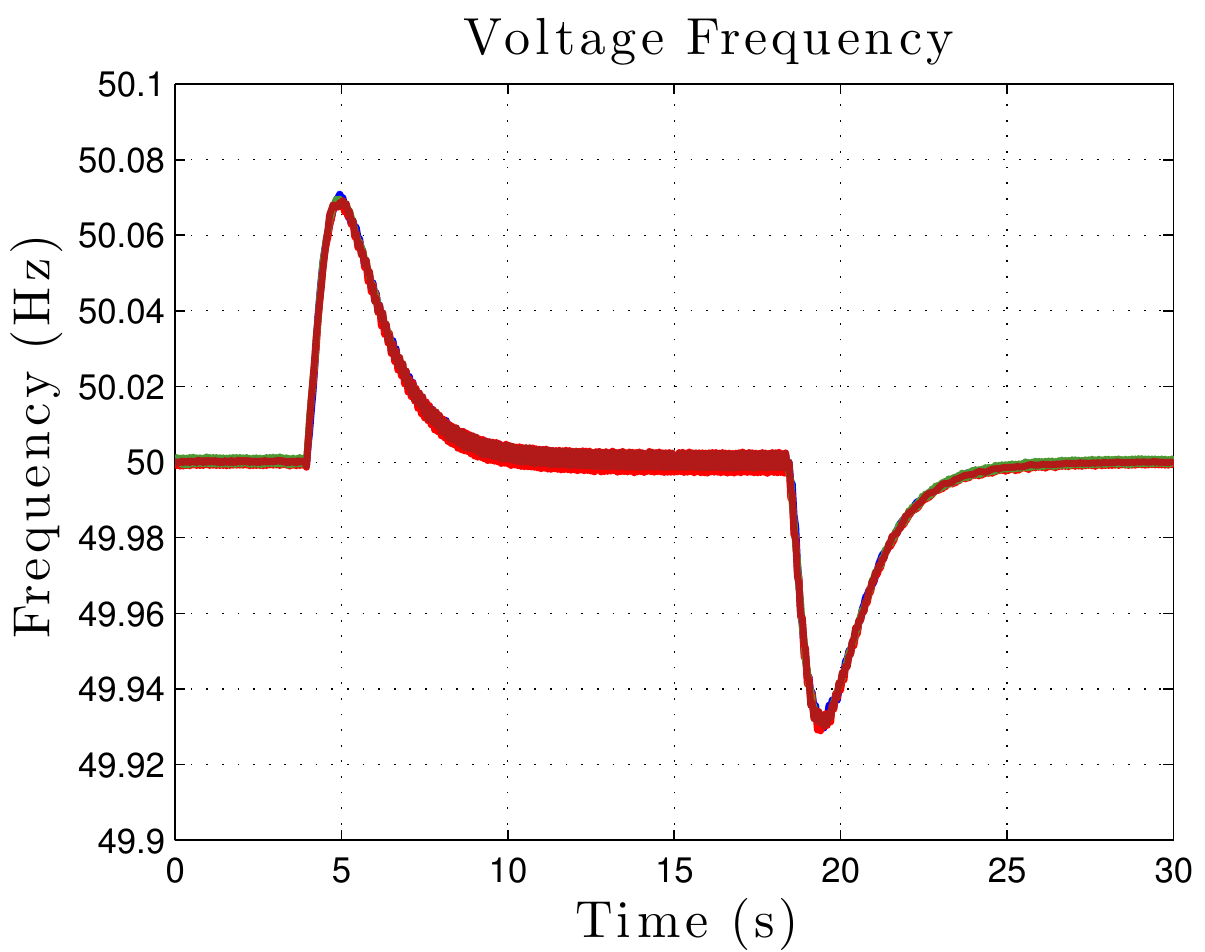}
                \label{Fig:5f}
        \end{subfigure}~
        \begin{subfigure}[!ht]{0.47\columnwidth}
                \includegraphics[width=\columnwidth]{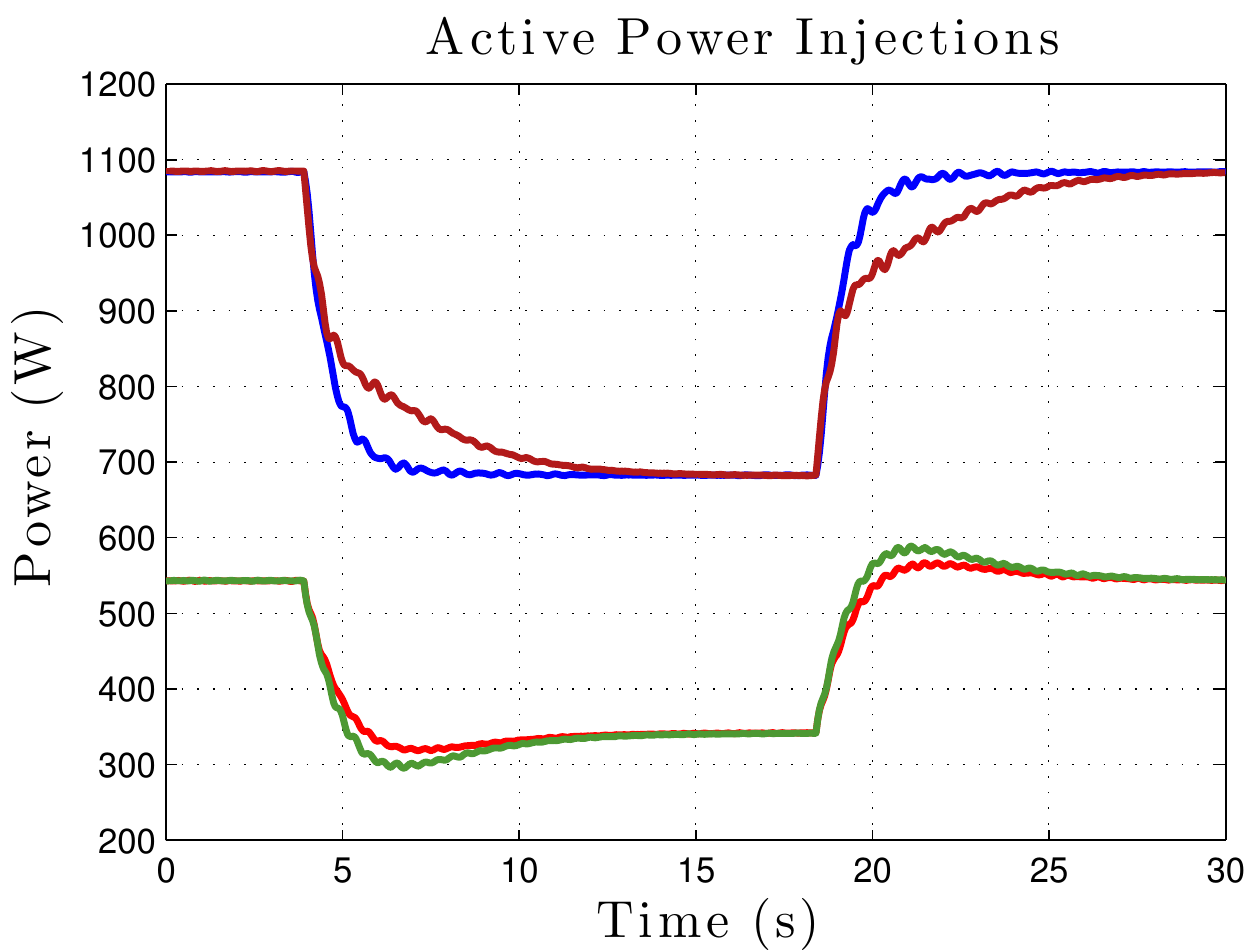}
                \label{Fig:5P}
                \end{subfigure}
%
		\captionsetup{justification=raggedright,singlelinecheck=false}
        \caption{Study 3 -- DAPI performance with heterogeneous controller gains. Control parameters are the same as in Study 1d.}\label{Fig:5}
\end{figure}

We examine the behavior of the frequency DAPI controller \eqref{Eq:Primary}--\eqref{Eq:Secondary} under inhomogeneous controller gains. Control parameters are the same as in Study 1d, except for variations in the integral gains $k_1 = 1.5$\,s, $k_2 = 1$\,s, $k_3 = 2$\,s and $k_4 = 0.5$\,s. The results are displayed in Figure \ref{Fig:5}. Note that the inhomogeneous controller gains leads to varying transient responses for the DGs, but the steady-state behavior and stability of the system is unchanged. This illustrates the utility of the gains $k_i$ and $\kappa_i$ in tuning the transient response of the DAPI-controlled microgrid.

\subsection{Study 4: Plug-and-Play Functionality}
\label{Sec:Study4}

\begin{figure}[!ht]
        \centering
        \begin{subfigure}[!ht]{0.47\columnwidth}
                \includegraphics[width=\columnwidth]{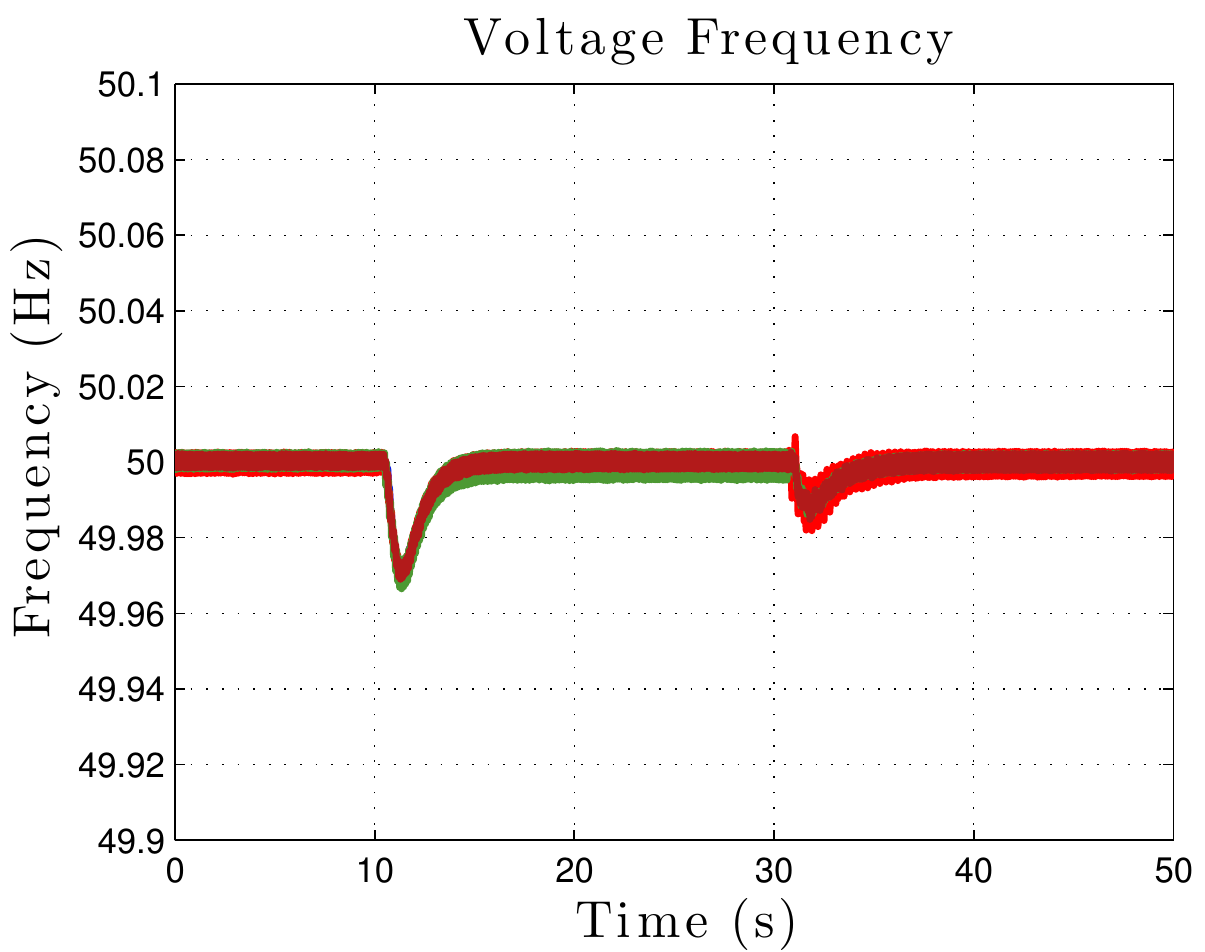}
                \label{Fig:4bf}
        \end{subfigure}~
        \begin{subfigure}[!ht]{0.47\columnwidth}
                \includegraphics[width=\columnwidth]{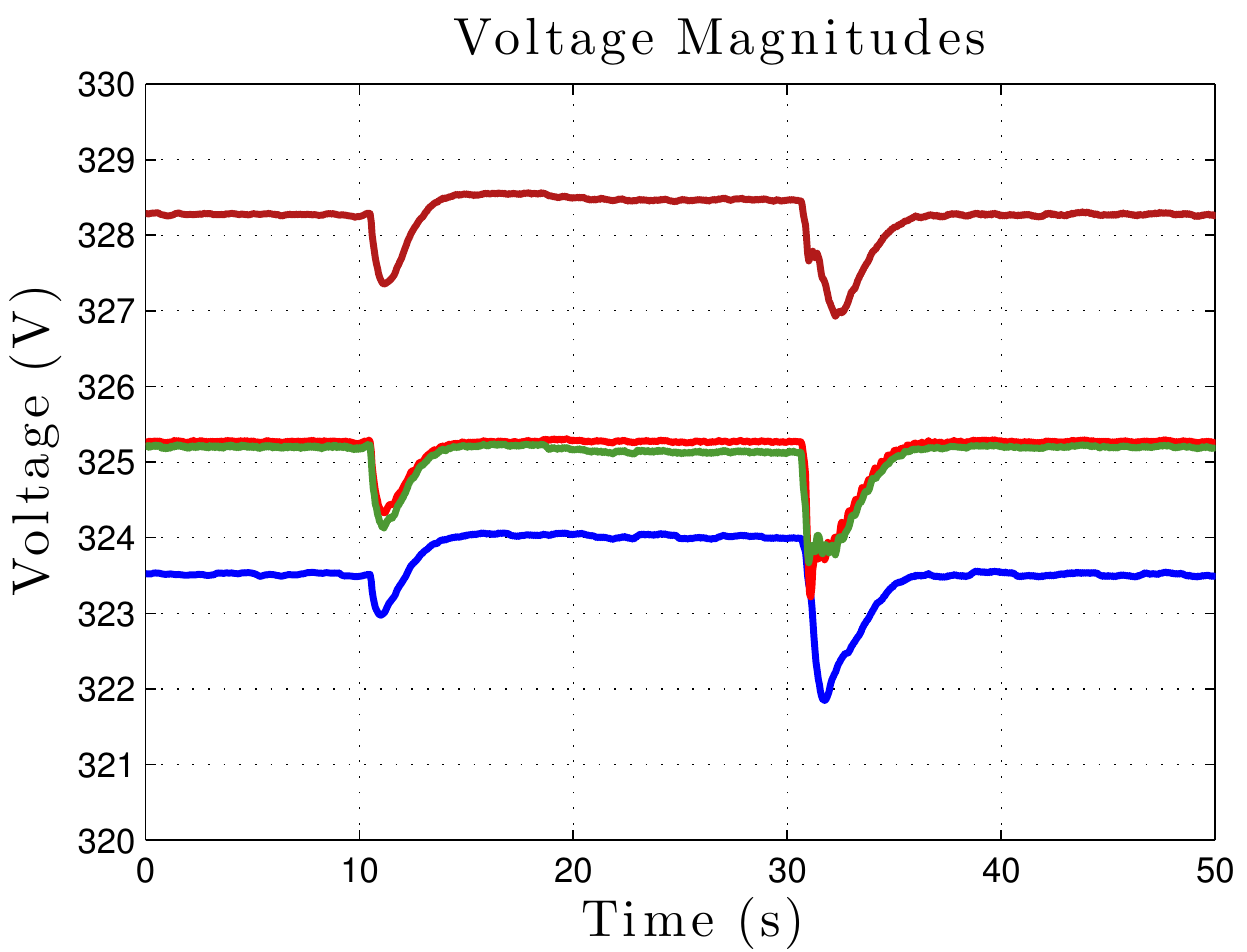}
                \label{Fig:4bE}
        \end{subfigure}
        
        \begin{subfigure}[!ht]{0.47\columnwidth}
                \includegraphics[width=\columnwidth]{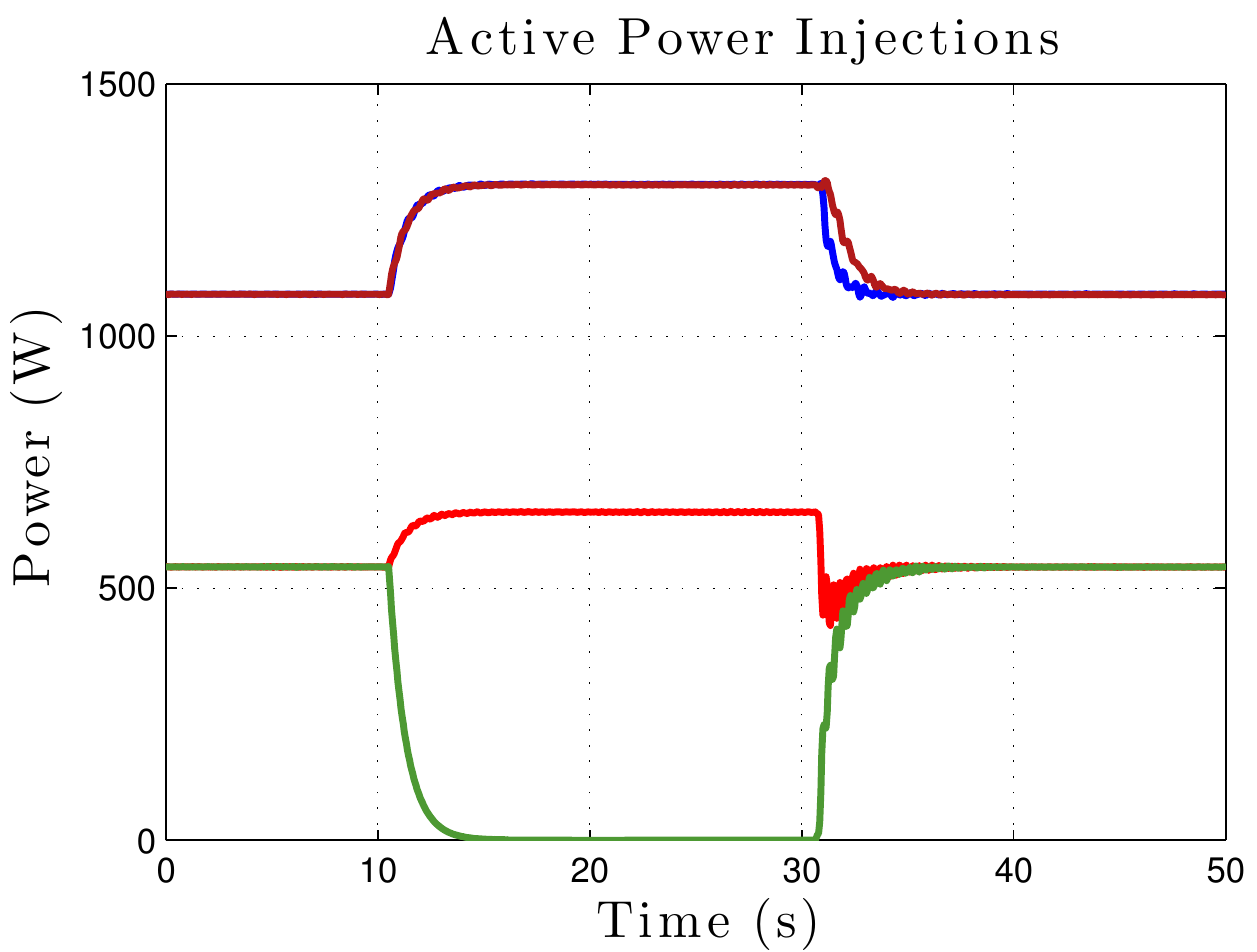}
                \label{Fig:4bP}
        \end{subfigure}~
        \begin{subfigure}[!ht]{0.47\columnwidth}
                \includegraphics[width=\columnwidth]{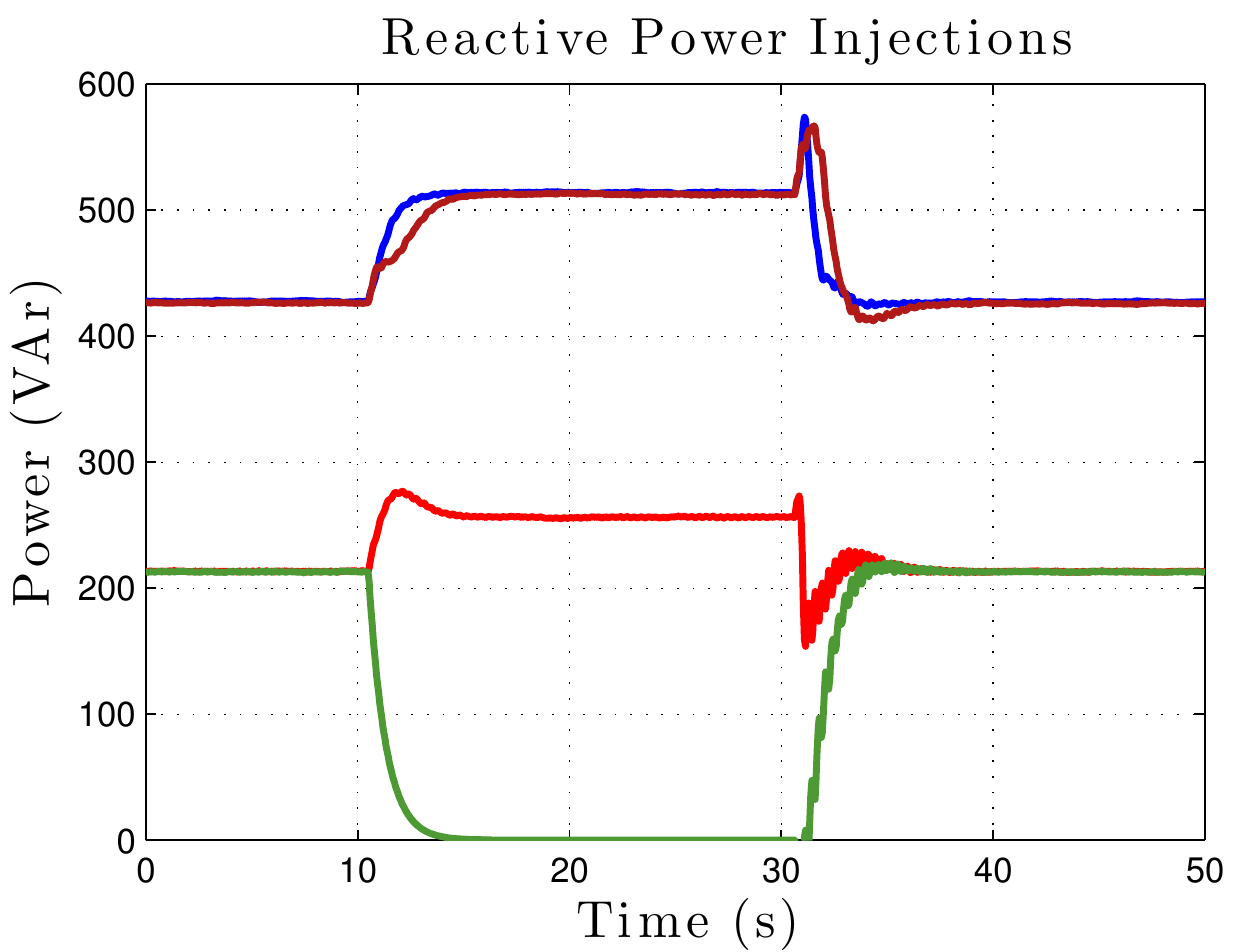}
                \label{Fig:4bQ}
        \end{subfigure}
%
		\captionsetup{justification=raggedright,singlelinecheck=false}
        \caption{Study 4 -- DAPI performance under plug-and-play operation. Control parameters are the same as in Study 1d.}\label{Fig:4b}
\end{figure}

The plug-and-play functionality of the controllers was tested by disconnecting unit 3 at $t = 10$\,s, and reconnecting it at $t = 30$\,s. A synchronization process was used in the downtime to synchronize unit 3 with the remaining microgrid before reconnection. Control parameters are the same as in Study 1d, and the results are displayed in Figure \ref{Fig:4b}. As in previous experiments, the DAPI controllers \eqref{Eq:Primary}--\eqref{Eq:SecondaryReactive} maintain accurate power sharing and frequency and voltage regulation before, during, and after the plug-and-play procedure, with minimal transients. The bus voltages and bus frequencies remain well regulated despite the disconnection of DG 3.


\section{Conclusions}~\label{Sec:Conclusions}~{ We have introduced a general distributed control methodology for primary/secondary control in islanded microgrids. By leveraging distributed averaging algorithms from multi-agent systems, the DAPI controllers achieve frequency regulation while sharing active power proportionally, and can be tuned to achieve either voltage regulation, reactive power sharing, or a compromise between the two. A small-signal stability analysis has been presented for the voltage DAPI controller along with a performance study, and the controllers have been validated through extensive experimental testing.

Large-signal stability of the voltage controller \eqref{Eq:PrimaryReactive}-\eqref{Eq:SecondaryReactive} remains an open analysis problem. Moreover, the secondary control goals of voltage regulation and reactive power sharing ignore an important factor for microgrid stability: the voltage levels at loads which are not collocated with DGs. 
An open problem is to design a controller guaranteeing that voltage levels at non-collocated load buses remain within tolerances while maintaining a high level of performance.}









\renewcommand{\baselinestretch}{1}
\bibliographystyle{IEEEtran}
\bibliography{brevalias_john,Main_Updates,New,Main,FB} 
%
%
%
%
%
%
\vspace{-2em}
\begin{IEEEbiography}[{\includegraphics[width=1in,height=1.25in,clip,keepaspectratio]{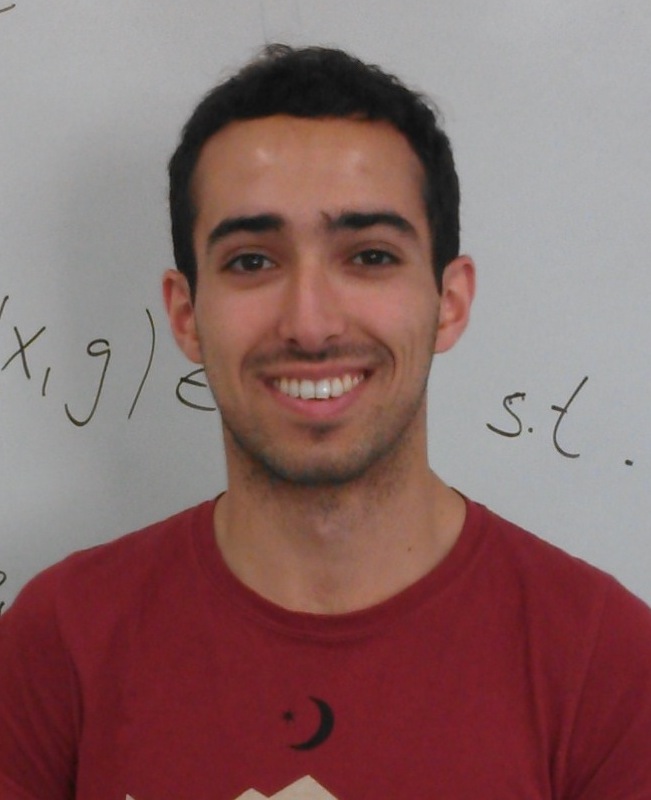}}]{John W. Simpson-Porco} (S'11) is a Ph.D. Candidate in the Department of Mechanical Engineering at the University of California Santa Barbara. He received his B.Sc. degree in Engineering Physics from Queen's University in 2010. Mr. Simpson-Porco is a recipient the 2012-2014 Automatica Best Paper Prize, the National Sciences and Engineering Research Council of Canada Fellowship, and the Center for Control, Dynamical Systems and Computation Outstanding Scholar Fellowship. His research interests are centered on the stability and control of multi-agent systems and complex dynamic networks, with a focus on modernized electric power grids.
\end{IEEEbiography}
\begin{IEEEbiography}[{\includegraphics[width=1in,height=1.25in,clip,keepaspectratio]{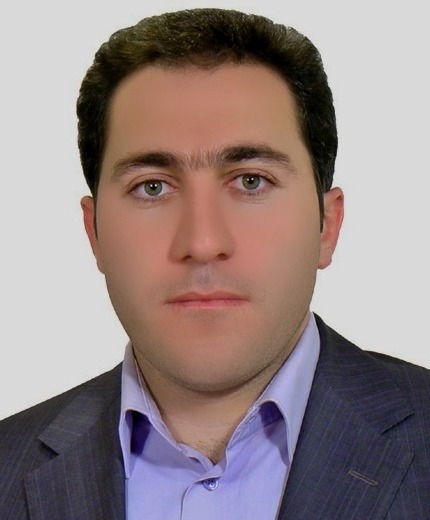}}]{Qobad Shafiee} (S'13--M'15) received the M.S. degree in electrical engineering from the Iran University of Science and Technology, Tehran, Iran, in 2007, and the PhD degree in power electronics and Microgrids from the Department of Energy Technology, Aalborg University, Aalborg, Denmark, in 2014. He worked with Department of Electrical and Computer Engineering, University of Kurdistan, Sanandaj, Iran, from 2007 to 2011, where he taught several electrical engineering courses and conducted research on load frequency control of power systems. From March 2014 to June 2014, he was a visiting researcher at the Electrical Engineering Department, University of Texas-Arlington, Arlington, TX, USA. Currently, he is a postdoctoral researcher at the Department of Energy Technology, Aalborg University. His main research interests include modeling, power management, hierarchical and distributed control applied to distributed generation in Microgrids.
\end{IEEEbiography}
%
\begin{IEEEbiography}[{\includegraphics[width=1in,height=1.25in,clip,keepaspectratio]{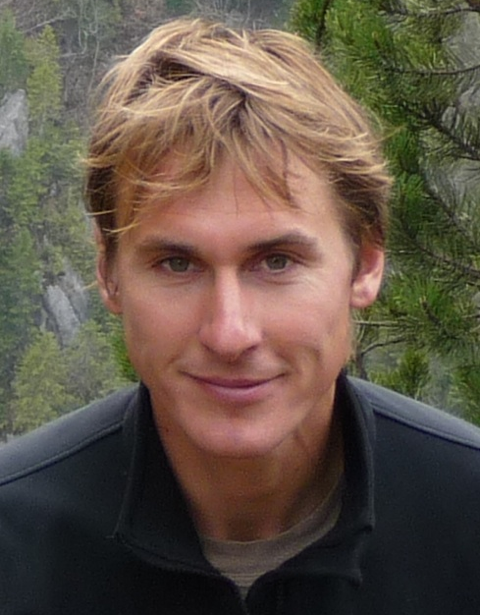}}]{Florian D\"{o}rfler} (S'09--M'13) is an Assistant Professor at the Automatic Control Laboratory at ETH ZŸrich. He received his Ph.D. degree in Mechanical Engineering from the University of California at Santa Barbara in 2013, and a Diplom degree in Engineering Cybernetics from the University of Stuttgart in 2008. From 2013 to 2014 he was an Assistant Professor at the University of California Los Angeles. His primary research interests are centered around distributed control, complex networks, and cyberÐphysical systems currently with applications in energy systems and smart grids. He is a recipient of the 2009 Regents Special International Fellowship, the 2011 Peter J. Frenkel Foundation Fellowship, the 2010 ACC Student Best Paper Award, the 2011 O. Hugo Schuck Best Paper Award, and the 2012-2014 Automatica Best Paper Award. As a co-advisor and a co-author, he has been a finalist for the ECC 2013 Best Student Paper Award.
\end{IEEEbiography}
%
%
%
%
\begin{IEEEbiography}[{\includegraphics[width=1in,height=1.25in,clip,keepaspectratio]{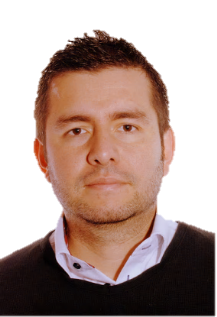}}]{Juan C. Vasquez} (M'12--SM'15) received the B.S. degree in electronics engineering from the Autonomous University of Manizales, Manizales, Colombia, and the Ph.D. degree in automatic control, robotics, and computer vision from the Technical University of Catalonia, Barcelona, Spain, in 2004 and 2009, respectively. He was with the Autonomous University of Manizales, where he has been teaching courses on digital circuits, servo systems and flexible manufacturing systems. In 2009, he worked as Post-doctoral Assistant with the Department of Automatic Control Systems and Computer engineering, Technical University of Catalonia, where he taught courses based on renewable energy systems, and power management on ac/dc minigrids and Microgrids.Ê Since 2011, he has been an Assistant Professor in Microgrids at the Department of Energy Technology, Aalborg University, Denmark, where he is co-responsible of the Microgrids research programme and co-advising more than 10 PhD students and a number of international visitors in research experience. He is a visiting scholar at the Center for Power Electronics Systems Ð CPES at Virginia Tech, Blacksburg, VA, USA. He has published more than 100 journal and conference papers and holds a pending patent. His current research interests include operation, energy management, hierarchical and cooperative control, energy management systems and optimization applied to Distributed Generation in AC/DC Microgrids. Dr Vasquez is member of the Technical Committee on Renewable Energy Systems TC-RES of the IEEE Industrial Electronics Society and the IEC System Evaluation Group SEG 4 work on LVDC Distribution and Safety for use in Developed and Developing Economies.
\end{IEEEbiography}

\begin{IEEEbiography}[{\includegraphics[width=1in,height=1.25in,clip,keepaspectratio]{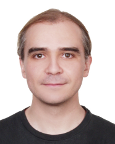}}]{Josep M. Guerrero} (S'01--M'04--SM'08--F'15) received the B.S. degree in telecommunications engineering, the M.S. degree in electronics engineering, and the Ph.D. degree in power electronics from the Technical University of Catalonia, Barcelona, in 1997, 2000 and 2003, respectively. Since 2011, he has been a Full Professor with the Department of Energy Technology, Aalborg University, Denmark, where he is responsible for the Microgrid Research Program. From 2012 he is a guest Professor at the Chinese Academy of Science and the Nanjing University of Aeronautics and Astronautics; from 2014 he is chair Professor in Shandong University; and from 2015 he is a distinguished guest Professor in Hunan University. His research interests is oriented to different microgrid aspects, including power electronics, distributed energy-storage systems, hierarchical and cooperative control, energy management systems, and optimization of microgrids and islanded minigrids. Prof. Guerrero is an Associate Editor for the IEEE TRANSACTIONS ON POWER ELECTRONICS, the IEEE TRANSACTIONS ON INDUSTRIAL ELECTRONICS, and the IEEE Industrial Electronics Magazine, and an Editor for the IEEE TRANSACTIONS on SMART GRID and IEEE TRANSACTIONS on ENERGY CONVERSION. He has been Guest Editor of the IEEE TRANSACTIONS ON POWER ELECTRONICS Special Issues: Power Electronics for Wind Energy Conversion and Power Electronics for Microgrids; the IEEE TRANSACTIONS ON INDUSTRIAL ELECTRONICS Special Sections: Uninterruptible Power Supplies systems, Renewable Energy Systems, Distributed Generation and Microgrids, and Industrial Applications and Implementation Issues of the Kalman Filter; and the IEEE TRANSACTIONS on SMART GRID Special Issue on Smart DC Distribution Systems. He was the chair of the Renewable Energy Systems Technical Committee of the IEEE Industrial Electronics Society. In 2014 he was awarded by Thomson Reuters as Highly Cited Researcher, and in 2015 he was elevated as IEEE Fellow for his contributions on Òdistributed power systems and microgrids.Ó
\end{IEEEbiography}

\begin{IEEEbiography}[{\includegraphics[width=1in,height=1.25in,clip,keepaspectratio]{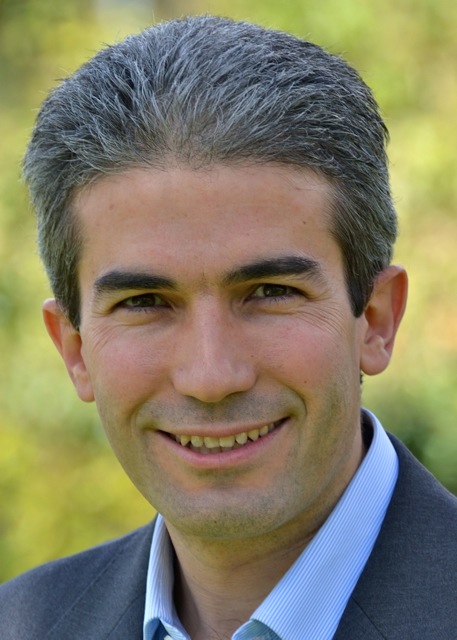}}]{Francesco Bullo} (S'95--M'99--SM'03Ð-F'10) is a Professor with the
Mechanical Engineering Department and the Center for Control, Dynamical
Systems and Computation at the University of California, Santa Barbara.  He
was previously associated with the University of Padova, the California
Institute of Technology and the University of Illinois at Urbana-Champaign.
His main research interests are network systems and distributed control
with application to robotic coordination, power grids and social networks.
He is the coauthor of "Geometric Control of Mechanical Systems" (Springer,
2004, 0-387-22195-6) and "Distributed Control of Robotic Networks"
(Princeton, 2009, 978-0-691-14195-4).  He received the 2008 IEEE CSM
Outstanding Paper Award, the 2010 Hugo Schuck Best Paper Award, the 2013
SIAG/CST Best Paper Prize, and the 2014 IFAC Automatica Best Paper
Award. He has served on the editorial boards of IEEE Transactions on
Automatic Control, ESAIM: Control, Optimization, and the Calculus of
Variations, SIAM Journal of Control and Optimization, and Mathematics of
Control, Signals, and Systems.
\end{IEEEbiography}




\end{document}